\documentclass[11pt,twoside, leqno]{article}

\usepackage{amssymb}
\usepackage{amsmath}

\allowdisplaybreaks

\def\rr{{\mathbb R}}
\def\rd{{{\rr}^d}}

\def\nn{{\mathbb N}}
\def\hh{{\mathbb H}}
\def\bbg{{\mathbb G}}

\def\cc{{\mathbb C}}

\def\cx{{\mathcal X}}

\def\cd{{\mathcal D}}
\def\ce{{\mathcal E}}

\def\cl{{\mathcal L}}

\def\cb{{\mathcal B}}

\def\fz{\infty}
\def\az{\alpha}
\def\supp{{\mathop\mathrm{\,supp\,}}}

\def\loc{{\mathop\mathrm{\,loc\,}}}
\def\lip{{\mathop\mathrm{\,Lip}}}
\def\lz{\lambda}
\def\dz{\delta}

\def\ez{\epsilon}

\def\kz{\kappa}
\def\bz{\beta}

\def\gz{{\gamma}}

\def\wz{\widetilde}

\def\hs{\hspace{0.3cm}}

\def\ls{\lesssim}
\def\gs{\gtrsim}

\def\hl{{H^1_\cl(\rn)}}

\def\lt{{L^2(\cx)}}

\def\lp{{L^p(\cx)}}

\def\hl{{\mathop\mathrm {HL}}}

\def\ati{{\mathrm {AOTI}}}
\def\bmo{{\mathop\mathrm{BMO}}}

\def\llip{{\lip_\cd}}
\def\blo{{\mathop\mathrm{BLO}}}

\def\esup{\mathop\mathrm{\,esssup\,}}
\def\einf{{\mathop{\mathrm{\,essinf\,}}}}
\def\dsum{\displaystyle\sum}
\def\diam{{\mathop\mathrm{\,diam\,}}}

\def\dint{\displaystyle\int}

\def\dfrac{\displaystyle\frac}
\def\dsup{\displaystyle\sup}
\def\dlim{\displaystyle\lim}

\def\r{\right}
\def\lf{\left}

\newtheorem{thm}{Theorem}[section]
\newtheorem{lem}{Lemma}[section]
\newtheorem{prop}{Proposition}[section]
\newtheorem{rem}{Remark}[section]
\newtheorem{cor}{Corollary}[section]
\newtheorem{defn}{Definition}[section]

\newtheorem{pf}{\bf Proof.}[section]

\numberwithin{equation}{section}

\begin{document}

\arraycolsep=1pt

\title{{\vspace{-4.5cm}\small\hfill\bf Nagoya Math. J., to appear}\\
\vspace{4cm}\Large\bf Localized Morrey-Campanato Spaces on Metric Measure
Spaces and Applications to Schr\"odinger Operators
\footnotetext{\hspace{-0.35cm} 2000 {\it
Mathematics Subject Classification}. {Primary 42B25; Secondary 42B35, 42B30.}
\endgraf{\it Key words and phrases.} Space of homogeneous type,
Heisenberg group, connected and simply connected nilpotent Lie group,
admissible function, Schr\"odinger operator,
atomic Hardy space, Morrey-Campanato space, maximal function, $g$-function.
\endgraf
 The first author is supported by the National
Natural Science Foundation (Grant No. 10871025) of China.}}
\author{Dachun Yang, Dongyong Yang and Yuan Zhou}
\date{ }
\maketitle

\begin{center}
\begin{minipage}{13.5cm}\small
{\noindent{\bf Abstract.} Let ${\mathcal X}$ be a space of
homogeneous type in the sense of Coifman and Weiss and ${\mathcal
D}$ a collection of balls in $\cx$. The authors introduce the
localized atomic Hardy space $H^{p,\,q}_{\mathcal D}({\mathcal X})$
with $p\in (0,1]$ and $q\in[1,\infty]\cap(p,\infty]$, the localized
Morrey-Campanato space ${\mathcal E}^{\alpha,\,p}_{\mathcal
D}({\mathcal X})$ and the localized Morrey-Campanato-BLO space
$\widetilde{\mathcal E}^{\alpha,\,p}_{\mathcal D}({\mathcal X})$
with $\az\in{\mathbb R}$ and $p\in(0, \infty)$ and establish their
basic properties including $H^{p,\,q}_{\mathcal D}({\mathcal
X})=H^{p,\,\infty}_{\mathcal D}({\mathcal X})$ and several
equivalent characterizations for ${\mathcal
E}^{\alpha,\,p}_{\mathcal D}({\mathcal X})$ and $\wz{\mathcal
E}^{\alpha,\,p}_{\mathcal D}({\mathcal X})$. Especially, the authors
prove that when $\az>0$ and $p\in[1,\,\fz)$, then
$\wz\ce^{\az,\,p}_{\mathcal D}(\mathcal X)=\ce^{\az,\,p}_{\mathcal
D}(\mathcal X)=\lip_{\mathcal D}(\alpha;\,\mathcal X)$, and when $p\in
(0,1]$, then the dual space of $H^{p,\,\infty}_{\mathcal
D}({\mathcal X})$ is ${\mathcal E}^{1/p-1,\,1}_{\mathcal
D}({\mathcal X})$. Let $\rho$ be an admissible function modeled on
the known auxiliary function determined by the Schr\"odinger
operator. Denote the spaces ${\mathcal E}^{\alpha,\,p}_{\mathcal
D}({\mathcal X})$ and $\widetilde{\mathcal E}^{\alpha,\,p}_{\mathcal
D}({\mathcal X})$, respectively, by ${\mathcal
E}^{\alpha,\,p}_{\rho}({\mathcal X})$ and $\widetilde{\mathcal
E}^{\alpha,\,p}_{\rho}({\mathcal X})$, when ${\mathcal D}$ is
determined by $\rho$. The authors then obtain the boundedness from
${\mathcal E}^{\alpha,\,p}_{\rho}({\mathcal X})$ to
$\widetilde{\mathcal E}^{\alpha,\,p}_{\rho}({\mathcal X})$ of the
radial and the Poisson semigroup maximal functions and the
Littlewood-Paley $g$-function which are defined via kernels modeled
on the semigroup generated by the Schr\"odinger operator. These
results apply in a wide range of settings, for instance, to the
Schr\"odinger operator or the degenerate Schr\"odinger operator on
${{\mathbb R}}^d$, or the sub-Laplace Schr\"odinger operator on
Heisenberg groups or connected and simply connected nilpotent Lie
groups. }
\end{minipage}
\end{center}

\section{Introduction\label{s1}}

\hskip\parindent The theory of Morrey-Campanato spaces plays an important
role in harmonic analysis and partial differential equations; see,
for example, \cite{c63, p69, tw, t92, s93, l07, n06, n08, dxy} and
their references. It is well-known that the dual space of
the Hardy space $H^p(\rd)$ with $p\in(0,\,1)$ is
the Morrey-Campanato space $\ce^{1/p-1, \,1}(\rd)$. Notice that Morrey-Campanato
spaces on $\rd$ are essentially related to the Laplacian $\Delta$,
where $\Delta\equiv\sum_{j=1}^d \frac{\partial^2}{\partial x_j^2}$.

On the other hand, there exists an increasing interest on the study of
Schr\"odinger operators on $\rd$ and the sub-Laplace Schr\"odinger
operators on connected and simply connected nilpotent Lie groups
with nonnegative potentials satisfying the reverse
H\"older inequality; see, for example,
\cite{f83,z99,s95,l99,dz99,dgmtz05,ll08,yz08,hl}. Let
$\cl\equiv -\Delta+V$ be the Schr\"odinger operator on $\rd$, where
the potential $V$ is a nonnegative locally integrable function.
Denote by $\cb_q(\rd)$ the class of functions satisfying the reverse
H\"older inequality of order $q$. For $V\in\cb_{d/2}(\rd)$ with
$d\ge3$, Dziuba\'nski et al \cite{dz99,dz03,dgmtz05} studied
the BMO-type space $\bmo_\cl(\rd)$ and the Hardy space $H^p_\cl(\rd)$
with $p\in(d/(d+1),1]$ and, especially, proved that the dual space
of $H^1_\cl(\rd)$ is $\bmo_\cl(\rd)$. Moreover,
they obtained the boundedness on these spaces of the variants of
several classical operators, including the radial maximal function
and the Littlewood-Paley $g$-function associated to $\cl$.
Recently, Huang and Liu \cite{hl}
further proved that the dual space of $H^p_\cl(\rd)$
is certain Morrey-Campanato space.
Let $\cx$ be an RD-space in \cite{hmy2}, which means that
$\cx$ is a space of homogeneous type in the sense of Coifman
and Weiss \cite{cw71, cw77} with the additional property
that a reverse doubling condition holds.
Let $\rho$ be a given admissible function
modeled on the known auxiliary function
determined by $V\in\cb_{d/2}(\rd)$
(see \cite{yz08} or \eqref{2.3} below).
Then the localized Hardy space $H^1_\rho(\cx)$,
the BMO-type space $\mathrm{\,BMO}_\rho({\mathcal X})$ and
the BLO-type space $\mathrm{\,BLO}_\rho({\mathcal X})$ were
introduced and studied by the authors of this paper in \cite{yz08,yyz}.
Moreover, the boundedness from $\mathrm{\,BMO}_\rho({\mathcal X})$
to $\mathrm{\,BLO}_\rho({\mathcal X})$ of
several maximal operators and the Littlewood-Paley $g$-function,
which are defined via kernels modeled on the
semigroup generated by the Schr\"odinger operator, was obtained
in \cite{yyz}.

The first purpose of this paper is to investigate behaviors
of these operators on localized Morrey-Campanato spaces on metric
measure spaces. To be precise, let ${\mathcal X}$ be a space of
homogeneous type, which is not necessary to be an RD-space,  and
${\mathcal D}$ be a collection of balls in $\cx$.
In Section \ref{s2} of this paper, we first introduce the
localized atomic Hardy space $H^{p,\,q}_{\mathcal D}({\mathcal X})$
with $p\in (0,1]$ and $q\in[1,\infty]\cap(p,\infty]$, the localized
Morrey-Campanato space ${\mathcal E}^{\alpha,\,p}_{\mathcal
D}({\mathcal X})$ and the localized Morrey-Campanato-BLO space
$\widetilde{\mathcal E}^{\alpha,\,p}_{\mathcal D}({\mathcal X})$
with $\az\in{\mathbb R}$ and $p\in(0, \infty)$, and establish their
basic properties including $H^{p,\,q}_{\mathcal D}({\mathcal
X})=H^{p,\,\infty}_{\mathcal D}({\mathcal X})$ and several
equivalent characterizations for ${\mathcal
E}^{\alpha,\,p}_{\mathcal D}({\mathcal X})$ and $\wz{\mathcal
E}^{\alpha,\,p}_{\mathcal D}({\mathcal X})$. Especially, we prove
that when $\az>0$ and $p\in[1,\,\fz)$, then
$\wz\ce^{\az,\,p}_\cd(\cx)=\ce^{\az,\,p}_\cd(\cx)=\lip_\cd(\az;\,\cx)$,
and when $p\in (0,1]$, then the dual space of $H^{p,\,\infty}_{\mathcal
D}({\mathcal X})$ is ${\mathcal E}^{1/p-1,\,1}_{\mathcal
D}({\mathcal X})$ (see Theorem \ref{t2.1} below).
Let $\rho$ be a given admissible function.
Modeled on the semigroup generated by the Schr\"odinger operator,
in Sections \ref{s3} and \ref{s4} of this paper,
we introduced the radial maximal operators $T^+$
and $P^+$ and Littlewood-Paley $g$-function $g(\cdot)$.
Then we establish the boundedness of $T^+$
and $P^+$ from $\ce_{\rho}^{\az,\,p}(\cx)$ to
$\wz\ce_{\rho}^{\az,\,p}(\cx)$ (see Theorems \ref{t3.1} and
\ref{t3.2} below). Here, for the set ${\mathcal D}$
determined by $\rho$, we denote
${\mathcal E}^{\alpha,\,p}_{\mathcal D}({\mathcal X})$
and $\widetilde{\mathcal E}^{\alpha,\,p}_{\mathcal D}({\mathcal
X})$, respectively, by ${\mathcal E}^{\alpha,\,p}_{\rho}({\mathcal
X})$ and $\widetilde{\mathcal E}^{\alpha,\,p}_{\rho}({\mathcal X})$.
Moreover, under the assumption that $g$-function $g(\cdot)$
is bounded on $L^p(\cx)$ with $p\in (1,\fz)$, we prove that
for every $f\in\ce_{\rho}^{\az,\,p}(\cx)$,
then $[g(f)]^2\in\wz\ce_{\rho}^{2\az,\,p/2}(\cx)$ with norm no more than
$C\|f\|_{\ce_{\rho}^{\az,\,p}(\cx)}^2$, where $C$ is a positive
constant independent of $f$ (see Theorem \ref{t4.1} below).
As a simple corollary of this, we obtain the
boundedness of $g(\cdot)$ from
$\ce_{\rho}^{\az,\,p}(\cx)$ to $\wz\ce_{\rho}^{\az,\,p}(\cx)$.
Notice that ${\mathcal E}^{0,\,p}_\rho({\mathcal X})=\bmo_\rho(\cx)$ and
$\wz{\mathcal E}^{0,\,p}_\rho({\mathcal X})=\blo_\rho(\cx)$
when $p\in[1,\fz)$. Thus, the results in this section
when $\az=0$ and $\cx$ is an RD-space were already obtained
in \cite{yyz}.

Finally, as the second purpose of this paper, in Section \ref{s5}
of this paper, we apply results obtained in Sections \ref{s3}
and \ref{s4} of this paper, respectively, to the Schr\"odinger operator or the
degenerate Schr\"odinger operator on $\rd$, the
sub-Laplace Schr\"odinger operator on Heisenberg groups or on
connected and simply connected nilpotent Lie
groups (see Propositions \ref{p5.1} through \ref{p5.5} below).
The nonnegative potentials of these Schr\"odinger operators are assumed to
satisfy the reverse H\"older inequality.

We now make some conventions. Throughout this
paper, we always use $C$ to denote a positive constant that is
independent of the main parameters involved but whose value may differ
from line to line. Constants with subscripts, such as $C_1$ and $A_1$, do not
change in different occurrences. If $f\le Cg$, we then write $f\ls
g$ or $g\gs f$; and if $f \ls g\ls f$, we then write $f\sim g.$
For any given ``normed" spaces $\mathcal A$  and $\mathcal B$,
the symbol ${\mathcal A}\subset {\mathcal B}$ means that for all $f\in \mathcal A$,
then $f\in\mathcal B$ and $\|f\|_{\mathcal B}\ls \|f\|_{\mathcal A}$.
We also use $B$ to denote a ball of $\cx$, and for $\lz>0$,
$\lz B$ denotes the ball with the same center as $B$, but radius
$\lz$ times the radius of $B$. Moreover, set $B^\complement\equiv\cx\setminus B$.
Also, for any set $E\subset\cx$, $\chi_E$ denotes its characteristic
function. For all $f\in L^1_\loc(\cx)$ and balls $B$, we always set
$f_{B}\equiv\frac1{\mu(B)}\int_Bf(y)\,d\mu(y)$.

\section{Localized Morrey-Campanato and Hardy spaces\label{s2}}

\hskip\parindent This section is divided into two subsections. In
Subsection \ref{s2.1}, we introduce the localized spaces
${\mathcal E}^{\alpha,\,p}_{\mathcal D}({\mathcal X})$ and
$\widetilde{\mathcal E}^{\alpha,\,p}_{\mathcal D}({\mathcal X})$
with $\az\in{\mathbb R}$ and $p\in(0, \infty)$,
we then establish the relations of these localized spaces
with their corresponding global versions and prove that
for all $\az\in[0, \fz)$ and $p\in (1, \fz)$,
$\ce^{\az,\,p}_\cd(\cx)=\ce^{\az,\,1}_\cd(\cx)$ and
$\wz\ce^{\az,\,p}_\cd(\cx)=\wz\ce^{\az,\,1}_\cd(\cx)$.
In Subsection \ref{s2.2}, we introduce the localized space
$H^{p,\,q}_{\mathcal D}({\mathcal X})$
with $p\in (0,1]$ and $q\in[1,\infty]\cap(p,\infty]$,
and show that $H^{p,\,q}_{\mathcal D}({\mathcal X})
=H^{p,\,\infty}_{\mathcal D}({\mathcal X})$
and the dual space of $H^{p,\,\infty}_{\mathcal D}({\mathcal X})$ is
${\mathcal E}^{1/p-1,\,1}_{\mathcal D}({\mathcal X})$.

\subsection{Localized Morrey-Campanato spaces\label{s2.1}}

\hskip\parindent We first recall the notion of spaces of
homogeneous type in \cite{cw71,cw77}.

\begin{defn}\label{d2.1} \rm
Let $(\cx,\, d)$ be a metric space endowed with a regular Borel measure
$\mu$ such that all balls defined by $d$ have finite and positive
measure. For any $x\in \cx $ and $r\in(0, \fz)$, set the ball
$B(x,r)\equiv\{y\in \cx :\ d(x,y)<r\}.$
The triple $(\cx,\,d,\,\mu)$ is called a
    space of homogeneous type if there exists
    a constant $A_1\in[1, \fz)$ such that for all $x\in \cx $ and $r\in(0, \fz)$,
    \begin{equation}\label{2.1}
    \mu(B(x, 2r))\le A_1\mu(B(x,r))\ (\mathrm{\it doubling\ property}).
    \end{equation}
\end{defn}

From \eqref{2.1}, it is not difficult to see that there exists positive constants
$A_2$ and $n$ such that for all
    $x\in \cx$, $r\in(0, \fz)$ and $\lz\in[1, \fz)$,
\begin{equation*}
 \mu(B(x,\lz r))\le A_2\lz^ n\mu(B(x,r)).
\end{equation*}

In what follows, we always set $V_r(x)\equiv\mu(B(x,\,r))$ and
$V(x,\,y)\equiv\mu(B(x,\,d(x,\,y)))$ for all $x,\,y\in\cx$ and
$r\in(0,\,\fz)$.

\begin{defn}\label{d2.2}\rm (\cite{yz08})
A positive function $\rho$ on $\cx$ is said to be admissible if
there exist positive constants $C_0$ and $k_0$ such that for all
$x,\,y\in\cx$,
\begin{equation}\label{2.2}
\frac1{\rho(x)} \le
C_0\frac1{\rho(y)}\lf(1+\frac{d(x,\,y)}{\rho(y)}\r)^{k_0}.
\end{equation}
\end{defn}

Obviously, if $\rho$ is a constant function, then $\rho$ is
admissible. Moreover, let $x_0\in\cx$ be fixed. The function
$\rho(y)\equiv (1+d(x_0,\,y))^s$ for all $y\in\cx$ with $s\in(-\fz,
1)$ also satisfies Definition \ref{d2.2} with $k_0=s/(1-s)$ when
$s\in[0, 1)$ and $k_0=-s$ when $s\in(-\fz, 0)$. Another non-trivial
class of admissible functions is given by the well-known reverse
H\"older class $\mathcal\cb_q(\cx, d, \mu)$, which is written
as $\mathcal\cb_q(\cx)$ for simplicity. Recall that a nonnegative potential $V$ is
said to be in $\cb_q(\cx)$ with $q\in(1,\,\fz]$ if there exists a
positive constant $C$ such that for all balls $B$ of $\cx$,
$$\lf(\frac1{|B|}\dint_B[V(y)]^q\,dy\r)^{1/q}
\le \frac{C}{|B|}\dint_BV(y)\,dy$$ with the usual modification made
when $q=\fz$. It is known that if $V\in \cb_q(\cx)$ for certain
$q\in(1,\,\fz]$, then $V$ is an $A_\fz(\cx)$ weight in the sense of
Muckenhoupt, and also $V\in \cb_{q+\ez}(\cx)$ for some $\ez\in(0,
\fz)$; see, for example, \cite{s93} and \cite{st89}. Thus
$\cb_q(\cx)=\cup_{q_1>q}\cb_{q_1}(\cx)$. For all $V\in \cb_q(\cx)$
with certain $q\in(1,\,\fz]$ and all $x\in\cx$, set
\begin{equation}\label{2.3}
\rho(x)\equiv[m(x, V)]^{-1}\equiv\sup\lf\{r>0:\hs \frac{r^2}
{\mu(B(x,\,r))}\dint_{B(x,\,r)}V(y)\,dy\le 1\r\};
\end{equation}
see, for example, \cite{s95} and also \cite{yz08}. It was also
proved in \cite{yz08} that $\rho$ in \eqref{2.3} is an admissible
function if $n\ge1$, $q>\max\{1,\,n/2\}$ and $V\in \cb_q(\cx)$.

We now recall the notion of Morrey-Campanato spaces and introduce
the definitions of Morrey-Campanato-BLO space and their localized versions.

\begin{defn}\label{d2.3}\rm
Let $\az\in\rr$ and $p\in(0,\,\fz)$.

(i)  A function $f\in L^p_{\loc}(\cx)$ is said to be in the Morrey-Campanato space
$\ce^{\az,\,p}(\cx)$ if
$$\|f\|_{\ce^{\az,\,p}(\cx)}\equiv\sup_{B\subset\cx}
\lf\{\frac1{[\mu(B)]^{1+p\az}}\int_B|f(y)-f_B|^pd\mu(y)\r\}^{1/p}<\fz,$$
where the supremum is taken over all balls $B\subset \cx$ and
$f_B=\frac1{\mu(B)}\int_Bf(z)\,d\mu(z)$.

(ii) A function $f\in L^p_{\loc}(\cx)$ is said to be in the Morrey-Campanato-BLO space
$\wz\ce^{\az,\,p}(\cx)$ if
$$\|f\|_{\wz\ce^{\az,\,p}(\cx)}\equiv
\sup_{B\subset \cx}\lf\{\frac1{[\mu(B)]^{1+p\az}}
\int_B\lf[f(y)-{\mathop\einf_B} f\r]^p\,d\mu(y)\r\}^{1/p}<\fz,$$
where the supremum is taken over all balls $B\subset \cx$.

(iii) Let $\az\in(0,\,\fz)$. A function $f$ on $\cx$ is said to be
in the Lipschitz space $\lip(\az;\,\cx)$ if there exists a nonnegative
constant $C$ such  that for all $x,\, y\in\cx$ and balls $B$
containing $x$ and $y$,
$$|f(x)-f(y)|\le C[\mu(B)]^\az.$$
The minimal nonnegative constant $C$ as above is called the norm
of $f$ in $\lip(\az;\,\cx)$ and denoted by
$\|f\|_{\lip(\az;\,\cx)}$.
\end{defn}

\begin{rem}\label{r2.1}\rm
(i) The space $\ce^{\az,\,p}(\cx)$ was first introduced by Campanato
in \cite{c63} when $\cx$ is a bounded subset of $\rd$ and $\mu$ is
the $d$-dimensional Lebesgue measure. When $\az=0$,
$\ce^{0,\,p}(\cx)$ is just the space $\bmo^p(\cx)$ (the space of functions of
bounded mean oscillation), and
$\ce^{0,\,p}(\cx)$ with $p\in[1,\fz)$ coincides with $\bmo^1(\cx)$;
see \cite{cw77}. For simplicity, we denote $\bmo^1(\cx)$ by $\bmo(\cx)$.

(ii) The space $\wz\ce^{0,\,p}(\cx)$ is just the space $\blo^p(\cx)$
(the space of functions of bounded lower oscillation). By
(i) of this remark and the fact that $\blo^1(\cx)\subset \bmo(\cx)$, it is easy
to see that $\wz\ce^{0,\,p}(\cx)$ with $p\in[1,\fz)$
coincides with $\blo^1(\cx)$. For simplicity, we denote
$\blo^1(\cx)$ by $\blo(\cx)$. Recall that $\blo(\cx)$ and
$\wz \ce^{\az,\,p}(\cx)$ are not linear spaces. The space
$\blo(\rd)$ was first introduced by Coifman and Rochberg \cite{cr80}
and $\wz\ce^{\az,\,p}(\rd)$ was introduced in \cite{hmy07}.

(iii) When $\az\in\rr$ and $p\in[1,\fz)$,
$\wz\ce^{\az,\,p}(\cx)\subset \ce^{\az,\,p}(\cx)$.
Moreover, when $\az\in(0,\,\fz)$ and $p\in[1,\fz)$, we have
$\wz\ce^{\az,\,p}(\cx)=\ce^{\az,\,p}(\cx)=\lip(\az;\,\cx)$ with
equivalent norms. In fact, Mac\'ias and Segovia
\cite{ms791} proved that when $\az\in(0, \fz)$ and
$p\in[1,\fz)$, $\ce^{\az,\,p}(\cx)=\lip(\az;\,\cx)$.
On the other hand, for all $f\in \ce^{\az,\,p}(\cx)$ and balls $B$,
$$\int_B [f(y)-{\mathop\einf_B}f]^p\,d\mu(y)
\le \int_B{\esup_{x\in B}}|f(y)-f(x)|^p\,d\mu(y)
\ls\|f\|^p_{\lip(\az;\,\cx)} [\mu(B)]^{1+p\az}, $$
which implies that $\|f\|_{\wz\ce^{\az,\,p}(\cx)}
\ls\|f\|_{\lip(\az;\,\cx)}\sim\|f\|_{\ce^{\az,\,p}(\cx)}.$
Thus, $\ce^{\az,\,p}(\cx)\subset\wz\ce^{\az,\,p}(\cx)$ and the claim
holds.
\end{rem}

\begin{defn}\label{d2.4}\rm
Let $\cd$ be a collection of balls in $\cx$, $p\in(0,\fz)$ and
$\az\in\rr$. Denote by $B$ any ball of $\cx$.

(i) A function $f\in L^p_{\loc}(\cx)$ is said to be in the localized Morrey-Campanato space
$\ce_\cd^{\az,\,p}(\cx)$ if
$$\begin{array}{cl}
\|f\|_{\ce_\cd^{\az,\,p}(\cx)}\equiv&\dsup_{B\notin\cd}
\lf\{\frac1{[\mu(B)]^{1+p\az}}\int_B|f(y)-f_B|^pd\mu(y)\r\}^{1/p}\\
&+\dsup_{B\in\cd}
\lf\{\frac1{[\mu(B)]^{1+p\az}}\int_B|f(y)|^pd\mu(y)\r\}^{1/p}<\fz,
\end{array}$$
where $f_B=\frac1{\mu(B)}\int_Bf(z)\,d\mu(z)$.

(ii) A function $f\in L^p_{\loc}(\cx)$ is said to be in the localized Morrey-Campanato-BLO space
$\wz\ce^{\az,\,p}_\cd(\cx)$ if
 \begin{eqnarray*}
\|f\|_{\wz\ce^{\az,\,p}_\cd(\cx)}&&\equiv \sup_{B\notin \cd}
\lf\{\frac1{[\mu(B)]^{1+p\az}}
\int_B\lf[f(y)-{\mathop\einf_B} f\r]^p\,d\mu(y)\r\}^{1/p}\\
&&\hs+\sup_{B\in\cd}\lf\{\frac1{[\mu(B)]^{1+p\az}}
\int_B|f(y)|^p\,d\mu(y)\r\}^{1/p}<\fz.
\end{eqnarray*}

(iii) Let $\az\in(0,\,\fz)$. A function $f$ on $\cx$ is said to be
in the localized Lipschitz space $\llip(\az;\,\cx)$ if there exists a nonnegative
constant $C$ such that for all $x,\, y\in\cx$ and balls $B$
containing $x$ and $y$ with $B\notin\cd$,
$$|f(x)-f(y)|\le C[\mu(B)]^\az,$$
and that for all balls $B\in\cd$, $\|f\|_{L^\fz(B)}\le
C[\mu(B)]^\az.$ The minimal nonnegative constant $C$ as above is
called the norm of $f$ in $\llip(\az;\,\cx)$ and denoted by
$\|f\|_{\llip(\az;\,\cx)}$.
\end{defn}

\begin{rem}\label{r2.2}\rm
(i) When $\az=0$ and $p\in[1,\fz)$, we denote $\ce_\cd^{0,\,p}(\cx)$
by $\bmo_\cd^p(\cx)$ and $\bmo_\cd^1(\cx)$ by $\bmo_\cd(\cx)$. And we
also denote  $\wz\ce_\cd^{0,\,p}(\cx)$ by $\blo^p_\cd(\cx)$
and $\wz\ce_\cd^{0,\,1}(\cx)$ by  $\blo_\cd(\cx)$. The localized
$\blo$ space was first introduced in \cite{hyy} in the setting of
$\rd$ endowed with a nondoubling measure.

(ii) If $\cx$ is the Euclidean space $\rd$ and
$\cd\equiv\{B(x,\,r):\ r\ge1,\,x\in\rd\}$, then $\bmo_\cd(\cx)$ is
just the localized ${\rm BMO}$ space of Goldberg \cite{g79}, and
$\lip_\cd(\az;\,\cx)$ with $\az\in(0,1)$ is just the inhomogeneous
Lipschitz space (see also \cite{g79}).

(iii) For all $\az\in\rr$ and $p\in (0,\fz)$,
$\wz\ce_\cd^{\az,\,p}(\cx)\subset \ce_\cd^{\az,\,p}(\cx)\subset
\ce^{\az,\,p}(\cx)$. For $\az\in(0,\,\fz)$,
$\lip_\cd(\az;\,\cx)\subset \lip(\az;\,\cx)$.

 (iv) Let $\rho$ be an admissible function
 and $\cd_\rho\equiv\{B(x,\,r):\ x\in\cx,\ r\ge \rho(x)\}$.
 In this case, we denote $\ce^{\az,\,p}_{\cd_\rho}(\cx)$, $\wz \ce^{\az,\,p}_{\cd_\rho}(\cx)$,
$\lip_{\cd_\rho}(\az;\,\cx)$, $\bmo_{\cd_\rho}(\cx)$ and $\blo_{\cd_\rho}(\cx)$,
respectively, by $\ce^{\az,\,p}_\rho(\cx)$, $\wz
\ce^{\az,\,p}_\rho(\cx)$, $\lip_\rho(\az;\,\cx)$, $\bmo_\rho(\cx)$
and $\blo_\rho(\cx)$. In \cite{yyz}, the spaces
$\bmo_\rho(\cx)$ and $\blo_\rho(\cx)$ when $\cx$
is an RD-space were introduced.
\end{rem}

The following results follow from Definitions \ref{d2.3} and
\ref{d2.4}.

\begin{lem}\label{l2.1} Let $\cd$ be a collection of balls in
$\cx$, $p\in[1,\fz)$ and $\az\in\rr$.

(i) Then $f\in \ce_\cd^{\az,\,p}(\cx)$ if and only if $f\in \ce
^{\az,\,p}(\cx)$ and $\sup_{B\in\cd}|f_B|[\mu(B)]^{-\az}<\fz$;
moreover,
\begin{equation*}
\|f\|_{\ce_\cd^{\az,\,p}(\cx)}\sim\|f\|_{\ce^{\az,\,p}(\cx)}+
\dsup_{B\in\cd}|f_B|[\mu(B)]^{-\az}.
\end{equation*}

(ii) Then $f\in \wz\ce_\cd^{\az,\,p}(\cx)$ if and only if $f\in \wz\ce
^{\az,\,p}(\cx)$ and $\sup_{B\in\cd}|f_B|[\mu(B)]^{-\az}<\fz$;
moreover,
\begin{equation*}
\|f\|_{\wz\ce_\cd^{\az,\,p}(\cx)}\sim\|f\|_{\wz\ce^{\az,\,p}(\cx)}+
\dsup_{B\in\cd}|f_B|[\mu(B)]^{-\az}.
\end{equation*}

 (iii) Let $\az\in(0,\fz)$. Then $f\in \llip(\az;\,\cx)$ if and only if
$f\in\lip(\az;\,\cx)$ and $\sup_{B\in\cd}[\mu(B)]^{-\az}$
$\|f\|_{L^\fz(B)}<\fz$ or
$\sup_{B\in\cd}|f_B|[\mu(B)]^{-\az}<\fz$; moreover,
\begin{equation*}
\begin{array}[t]{cl}
\|f\|_{\llip(\az;\,\cx)}&\sim \|f\|_{\lip(\az;\,\cx)}+
\dsup_{B\in\cd}\|f\|_{L^\fz(B)}[\mu(B)]^{-\az}\\
&\sim\|f\|_{\lip(\az;\,\cx)}+\dsup_{B\in\cd}|f_B|[\mu(B)]^{-\az}.
\end{array}
\end{equation*}
\end{lem}

\begin{pf}\rm
We first prove (i). If $f\in \ce ^{\az,\,p}(\cx)$ and
$\sup_{B\in\cd}|f_B|[\mu(B)]^{-\az}<\fz$, from Definitions
\ref{d2.3} and \ref{d2.4}, it follows that
\begin{equation}\label{2.4}
\|f\|_{\ce_\cd^{\az,\,p}(\cx)}\le2\|f\|_{\ce^{\az,\,p}(\cx)}
+\dsup_{B\in\cd}|f_B|[\mu(B)]^{-\az}.
\end{equation}
Conversely, if $f\in \ce_\cd^{\az,\,p}(\cx)$, then by the H\"older
inequality, we have
\begin{eqnarray*}
\|f\|_{\ce^{\az,\,p}(\cx)}+\dsup_{B\in\cd}|f_B|[\mu(B)]^{-\az} &&\le
\|f\|_{\ce_\cd^{\az,\,p}(\cx)}+2\dsup_{B\in\cd}|f_B|[\mu(B)]^{-\az}
\le3\|f\|_{\ce_\cd^{\az,\,p}(\cx)},
\end{eqnarray*}
which together with \eqref{2.4} gives (i).

The proofs of (ii) and (iii) are similar. We omit the details, which
completes the proof of Lemma \ref{l2.1}.
\end{pf}

\begin{lem}\label{l2.2} Let $\cd$ be a collection of balls in
$\cx$ and $p\in[1,\fz)$.

(i) Then $\bmo_\cd(\cx)=\bmo_\cd^p(\cx)$ and
$\blo_\cd(\cx)=\blo_\cd^p(\cx)$ with equivalent norms.

(ii) When $\az\in(0,\,\fz)$, $\wz\ce_\cd^{\az,\,p}(\cx)=\ce_\cd^{\az,\,p}(\cx)
=\lip_\cd(\az;\,\cx)$ with equivalent norms.
\end{lem}

\begin{pf}\rm To prove (i), we first assume that $f\in\bmo_\cd^p(\cx)$.
Then by the H\"older
inequality, we have $f\in\bmo_\cd(\cx)$ and
$\|f\|_{\bmo_\cd(\cx)}\le\|f\|_{\bmo_\cd^p(\cx)}$. Conversely, if
$f\in\bmo_\cd(\cx)$, then from Lemma \ref{l2.1} (i) with $\az=0$,
Remark \ref{r2.1} (i) and Remark \ref{r2.2} (iii), it follows that
\begin{eqnarray*}
\|f\|_{\bmo_\cd^p(\cx)}\ls\|f\|_{\bmo^p(\cx)}+ \sup_{B\in\cd}|f_B|
&&\ls\|f\|_{\bmo_\cd(\cx)},
\end{eqnarray*}
which implies that $f\in\bmo_\cd^p(\cx)$ and
$\|f\|_{\bmo_\cd^p(\cx)}\ls\|f\|_{\bmo_\cd(\cx)}$. Thus
$\bmo_\cd(\cx)=\bmo_\cd^p(\cx)$ with equivalent norms.
The proof for $\blo_\cd(\cx)=\blo_\cd^p(\cx)$
is similar and we omit the details.

To prove (ii), by  Lemma \ref{l2.1} and Remark \ref{r2.1} (iii), we obtain
\begin{eqnarray*}
\|f\|_{\ce_\cd^{\az,\,p}(\cx)}&&\sim
\|f\|_{\ce^{\az,\,p}(\cx)}+\dsup_{B\in\cd}|f_B|[\mu(B)]^{-\az}\\
&&\sim
\|f\|_{\wz\ce^{\az,\,p}(\cx)}+\dsup_{B\in\cd}|f_B|[\mu(B)]^{-\az}
\sim\|f\|_{\wz\ce_\cd^{\az,\,p}(\cx)}\\
&&\sim\|f\|_{\lip(\az;\,\cx)}
+\dsup_{B\in\cd}|f_B|[\mu(B)]^{-\az}\sim\|f\|_{\llip(\az;\,\cx)},
\end{eqnarray*}
which implies (ii). This finishes the proof of Lemma \ref{l2.2}.
\end{pf}

The space $\cx$ is said to have the reverse doubling property if
there exist constants $\kz\in(0, n]$ and $A_3\in(0, 1]$ such that
for all $x\in \cx$, $r\in(0, \diam(\cx)/2]$ and $\lz\in[1,
\diam(\cx)/(2r)]$,
\begin{equation}\label{2.5}
A_3\lz^\kz\mu(B(x,r))\le\mu(B(x,\lz r)).
\end{equation}
If $(\cx,d,\mu)$ satisfies the conditions \eqref{2.1} and \eqref{2.5},
then $(\cx,d,\mu)$ is called an RD-space (see \cite{hmy2}).

\begin{lem}\label{l2.3}
Let $\cx$ be an RD-space, $\rho$ an admissible
function on $\cx$ and $\cd_\rho$ as in Remark \ref{r2.2} (iv).
If $\az\in(-\fz,\, 0)$ and $p\in[1,\,\fz)$, then
\begin{eqnarray*}
\|f\|_{\ce^{\az,\,p}_{\rho}(\cx)} &&\sim \sup_{B\subset
\cx}\lf\{\frac1{[\mu(B)]^{1+\az p}}\dint_B|f(x)|^p\,d\mu(x)\r\}^{1/p}.
\end{eqnarray*}
\end{lem}

\begin{pf}\rm
 An application of the H\"older inequality leads
to that
$$\|f\|_{\ce^{\az,\,p}_{\rho}(\cx)}\ls \sup_{B\subset \cx}
\lf\{\frac1{[\mu(B)]^{1+\az p}}\dint_B|f(x)|^p\,d\mu(x)\r\}^{1/p}.$$
Conversely, if $B\in{\cd_\rho}$, then by Definition \ref{d2.4} (i),
we have
$$\lf\{\frac1{[\mu(B)]^{1+\az p}}\dint_B|f(x)|^p\,d\mu(x)\r\}^{1/p}
\le \|f\|_{\ce^{\az,\,p}_{\rho}(\cx)}.$$ Now we assume that $B\equiv
B(x_0, r)\notin {\cd_\rho}$. Let $J_0\in \nn$ such that
$2^{J_0-1}r<\rho(x_0)\le 2^{J_0}r$. From $\az\in(-\fz, 0)$,
\eqref{2.1}, \eqref{2.5} and the H\"older inequality, it follows that
\begin{eqnarray*}
&&\lf\{\frac1{[\mu(B)]^{1+\az p}}\dint_B|f(x)|^p\,d\mu(x)\r\}^{1/p}\\
&&\hs\le\frac1{[\mu(B)]^\az}\lf\{\lf[\frac1{\mu(B)}\dint_B|f(x)-f_B|^p\,d\mu(x)\r]^{1/p}
+\sum_{j=1}^{J_0}\lf|f_{2^{j-1}B}-f_{2^jB}\r|+\lf|f_{2^{J_0}B}\r|\r\}\\
&&\hs\ls\lf(1+\sum_{j=1}^{J_0}2^{j\kappa\az}\r)\|f\|_{\ce^{\az,\,p}_{\rho}(\cx)}
\ls\|f\|_{\ce^{\az,\,p}_{\rho}(\cx)},
\end{eqnarray*}
which completes the proof of Lemma \ref{l2.3}.
\end{pf}

\begin{rem}\label{r2.3}\rm
Let $\cx$ be an RD-space and $p\in [1,\,\fz)$.

(i) Then Lemma \ref{l2.3}
implies that $\ce^{\az,\,p}_\rho(\cx)$ with
$\az\in(-1/p, 0)$ coincides with the so-called
Morrey space (see, for example, \cite{p69,t92}
for the case $\cx=\rd$).

(ii) Let $\az<0$. For all $f\ge0$, $f\in\ce_\cd^{\az,\,p}(\cx)$
if and only if $f\in \wz\ce_\cd^{\az,\,p}(\cx)$ and
moreover, $\|f\|_{\wz\ce_\cd^{\az,\,p}(\cx)}\sim
\|f\|_{\ce_\cd^{\az,\,p}(\cx)}$.
In fact, by Remark \ref{r2.2} (iii), we only need to show
that for all $f\ge0$, $f\in\ce_\cd^{\az,\,p}(\cx)$
implies that $f\in \wz\ce_\cd^{\az,\,p}(\cx)$ and
$\|f\|_{\wz\ce_\cd^{\az,\,p}(\cx)}\ls
\|f\|_{\ce_\cd^{\az,\,p}(\cx)}$.  By Lemma 2.3, the fact that
$\az<0$ and $f\ge0$, we see that for all balls $B\notin\cd$,
\begin{eqnarray*}
\int_B\lf[f(x)-\mathop\einf_{B}f\r]^p\,d\mu(x)\le
\int_B\lf[f(x)\r]^p\,d\mu(x)\ls[\mu(B)]^{1+\az p}\|f\|^p_{\ce_\cd^{\az,\,p}(\cx)},
\end{eqnarray*}
which implies the claim.

(iii) If $\cx$ is not an RD-space, it is not clear if
Lemma \ref{l2.3} still holds.
\end{rem}

We also have the following conclusions which are
used in Sections \ref{s3} and \ref{s4}.

\begin{lem}\label{l2.4}
Let $\az\in\rr$, $p\in[1,\fz)$, $\rho$ be an admissible function
on $\cx$ and $\cd_\rho$ as in Remark \ref{r2.2} (iv).
Then there exists a positive constant $C$ such that for
all $f\in \ce^{\az,\,p}_{\rho}(\cx)$,

(i) for all balls $B\equiv B(x_0,\,r)\not\in {\cd_\rho}$,
$$\frac1{\mu(B)}\int_B|f(z)|\,d\mu(z)\le\lf\{
\begin{array}{ll}
 C\lf(\frac{\rho(x_0)}{r}\r)^{\az n }
[\mu(B)]^\az\|f\|_{\ce^{\az,\,p}_{\rho}(\cx)}, \hs &\az>0,\\
C\lf(1+\log\frac{\rho(x_0)}{r}\r)[\mu(B)]^\az
\|f\|_{\ce^{\az,\,p}_{\rho}(\cx)},\hs &\az\le 0;
\end{array}\r.$$

(ii) for all $x\in\cx$ and $0<r_1<r_2$,

$$|f_{B(x,\,r_1)}-f_{B(x,\,r_2)}|\le
\lf\{\begin{array}{ll} C \lf(\frac{r_2}{r_1}\r)^{\az n }
[\mu(B(x,\,r_1)]^\az\|f\|_{\ce^{\az,\,p}_{\rho}(\cx)},&\hs \az>0,\\
C\lf(1+\log\frac{r_2}{r_1}\r)[\mu(B(x,\,r_1)]^\az\|f\|_{\ce^{\az,\,p}_{\rho}(\cx)},&\hs\az\le0.
\end{array}\r.$$
\end{lem}

\begin{pf}\rm
If (ii) holds, then by the H\"older inequality, we see that
for all $f\in \ce^{\az,\,p}_{\rho}(\cx)$
and $B\notin {\cd_\rho}$,
\begin{eqnarray*}
\frac1{\mu(B)}\int_B|f(x)|\,d\mu(x)&&\le
\frac1{\mu(B)}\int_B|f(x)-f_B|\,d\mu(x)+
\lf|f_B-f_{B(x_0,\,\rho(x_0))}\r|\\
&&\hs+
\frac1{\mu(B(x_0,\,\rho(x_0)))}\int_{B(x_0,\,\rho(x_0))}|f(x)|\,d\mu(x)\\
&&\ls\lf\{[\mu(B)]^\az+ [\mu(B(x_0, \rho(x_0)))]^\az\r\}
\|f\|_{\ce^{\az,\,p}_{\rho}(\cx)}
+\lf|f_B-f_{B(x_0,\,\rho(x_0))}\r|.
\end{eqnarray*}
Then (i) follows from this fact together with \eqref{2.1},
$r<\rho(x_0)$
(because $B\notin {\cd_\rho}$) and (ii).

To prove (ii), let $j_0$ be the smallest integer such that
$2^{j_0}r_1\ge r_2$. Another application of \eqref{2.1} leads to
that
\begin{eqnarray*}
\lf|f_{B(x,\,2^{j_0}r_1)}-f_{B(x,\,r_2)}\r|
&&\ls \frac1{\mu(B(x,\,2^{j_0}r_1))}
\int_{B(x,\,2^{j_0}r_1)}\lf|f-f_{B(x,\,2^{j_0}r_1)}\r|\,d\mu(z)\\
&&\ls \lf[\mu(B(x,\,2^{j_0}r_1))\r]^\az
\|f\|_{\ce^{\az,\,p}_{\rho}(\cx)}.
\end{eqnarray*}
Similarly, we see that for all $j\in\nn\cup\{0\}$,
$$\lf|f_{B(x,\,2^{j}r_1)}-f_{B(x,\,2^{j+1}r_1)}\r| \ls
\lf[\mu\lf(B\lf(x,\,2^{j+1}r_1\r)\r)\r]^\az
\|f\|_{\ce^{\az,\,p}_{\rho}(\cx)}.$$ Then we have
\begin{eqnarray*}
\lf|f_{B(x,\,r_1)}-f_{B(x,\,r_2)}\r|&&\ls\sum_{j=0}^{j_0-1}
\lf|f_{B(x,\,2^{j}r_1)}-f_{B(x,\,2^{j+1}r_1)}\r|
+\lf|f_{B(x,\,2^{j_0}r_1)}-f_{B(x,\,r_2)}\r|\\
&&\ls\sum_{j=0}^{j_0-1}[\mu(B(x,\,2^{j+1}r_1))]^\az
\|f\|_{\ce^{\az,\,p}_{\rho}(\cx)}.
\end{eqnarray*}
If $\az\in(-\fz, 0]$, from the choice of $j_0$, we deduce that
$$|f_{B(x,\,r_1)}-f_{B(x,\,r_2)}|
\ls\lf(1+\log\frac{r_2}{r_1}\r)[\mu(B(x,\,r_1))]^\az\|f\|_{\ce^{\az,\,p}_{\rho}(\cx)};$$
if $\az\in(0, \fz)$, by \eqref{2.1}, we obtain that
$$|f_{B(x,\,r_1)}-f_{B(x,\,r_2)}|
\ls\lf(\frac{r_2}{r_1}\r)^{\az n }[\mu(B(x,\,r_1))]^\az
\|f\|_{\ce^{\az,\,p}_{\rho}(\cx)}.$$
This finishes the proof of Lemma \ref{l2.4}.
\end{pf}

\subsection{Localized Hardy spaces\label{s2.2}}

\hskip\parindent We begin with the notion of atoms.

\begin{defn}\label{d2.5}\rm
Let $\cd$ be a collection of balls in $\cx$, $p\in(0,1]$ and $q\in
[1,\fz]\cap(p,\fz]$.

(i) A function $a$ supported in a ball $B\subset\cx$ is called a
$(p,\,q)$-atom if $\int_\cx a(x)\,d\mu(x)=0$ and
 $\|a\|_{L^q(\cx)}\le[\mu(B)]^{1/q-1/p}$ (see \cite{cw77}).

(ii) A function $b$ supported in a ball $B\in\cd$ is called a
$(p,\,q)_\cd$-atom if $\|b\|_{L^q(\cx)}\le[\mu(B)]^{1/q-1/p}.$
\end{defn}

\begin{rem}\label{r2.4}\rm
(i) Every $(1,\,q)$-atom or $(1,\,q)_\cd$-atom $a$ belongs to
$L^1(\cx)$ with $\|a\|_{L^1(\cx)}\le1$.

(ii) Let $p\in(0,1)$. If $a$ is a $(p,\,q)$-atom, then $a\in
(\lip(1/p-1;\,\cx))^\ast\subset (\llip(1/p-1;\,\cx))^\ast$ and
$\|a\|_{(\llip(1/p-1;\,\cx))^\ast}\le\|a\|_{(\lip(1/p-1;\,\cx))^\ast}\le1$;
if $b$ is a $(p,\,q)_\cd$-atom, then $b\in
(\llip(1/p-1;\,\cx))^\ast$ and
$\|b\|_{(\llip(1/p-1;\,\cx))^\ast}\le1$.

\end{rem}

\begin{defn}\label{d2.6}\rm(\cite{cw77})
Let $p\in(0,1]$ and $q\in [1,\fz]\cap(p,\fz]$. A function $f\in
L^1(\cx)$ or a linear functional $f\in(\lip(1/p-1;\,\cx))^\ast$ when
$p\in(0,1)$ is said to be in the Hardy space $H^{1,\,q}(\cx)$ when $p=1$
or in $H^{p,\,q}(\cx)$ when $p\in(0,1)$ if there exist $(p,\,q)$-atoms
$\{a_j\}_{j=1}^\fz$ and $\{\lambda_j\}_{j=1}^\fz\subset\cc$ such
that $f=\sum_{j\in\nn}\lambda_ja_j,$ which converges in $L^1(\cx)$
when $p=1$ or in $(\lip(1/p-1;\,\cx))^\ast$ when $p\in(0,1)$,
and $\sum_{j\in\nn}|\lambda_j|^p<\fz$. Moreover, the norm of $f$ in
$H^{p,\,q}(\cx)$ with $p\in(0,1]$ is defined by
$$\|f\|_{H^{p,\,q}(\cx)}
\equiv\inf\Bigg\{\bigg(\sum_{j\in\nn}|\lambda_j|^p\bigg)^{1/p}\Bigg\},$$
where the infimum is taken over all the above decompositions of $f$.
\end{defn}

\begin{rem}\label{r2.5}\rm
Coifman and Weiss \cite{cw77} proved that $H^{p,\,q}(\cx)$ and
$H^{p,\,\fz}(\cx)$ coincide with equivalent norms for all
$p\in(0,1]$ and $q\in[1,\,\fz)\cap(p,\,\fz)$. Thus, for all
$p$, $q$ in this range, we denote
$H^{p,\,q}(\cx)$ simply by $H^p(\cx)$. We remark that
Coifman and Weiss \cite{cw77} also proved that the dual space of $H^p(\cx)$ is
$\bmo(\cx)$ when $p=1$ or $\lip(1/p-1;\,\cx)$ when $p\in(0, 1)$.
\end{rem}

\begin{defn}\label{d2.7}\rm
Let $\cd$ be a collection of balls in $\cx$, $p\in(0,1]$ and $q\in
[1,\fz]\cap(p,\fz]$. A function $f\in L^1(\cx)$ or a linear
functional $f\in(\llip(1/p-1;\,\cx))^\ast$ when $p\in(0,1)$ is said
to be in $H_\cd^{1,\,q}(\cx)$ when $p=1$ or
$H_\cd^{p,\,q}(\cx)$ when $p\in(0,1)$ if there exist
$\{\lambda_j\}_{j\in\nn}$, $\{\nu_k\}_{k\in\nn}\subset\cc$,
$(p,\,q)$-atoms $\{a_j\}_{j\in\nn}$ and $(p,\,q)_\cd$-atoms
$\{b_k\}_{k\in\nn}$
 such that
$$f=\sum_{j\in\nn}\lambda_ja_j+
\sum_{k\in\nn}\nu_kb_k,$$ which converges in $L^1(\cx)$ when $p=1$
or in $(\llip(1/p-1;\,\cx))^\ast$ when $p\in(0, 1]$, and
$\sum_{j\in\nn}|\lambda_j|^p+ \sum_{k=1}^\fz|\nu_k|^p<\fz$.
Moreover, the norm of $f$ in $H_\cd^{p,\,q}(\cx)$ is defined by
$$\|f\|_{H_\cd^{p,\,q}(\cx)}
\equiv\inf\Bigg\{\bigg(\sum_{j\in\nn}|\lambda_j|^p+
\sum_{k\in\nn}|\nu_k|^p\bigg)^{1/p}\Bigg\},$$ where the infimum is
taken over all the above decompositions of $f$.
\end{defn}

\begin{rem}\label{r2.6}\rm
Let $p\in(0, 1]$ and $q\in[1,\fz]\cap(p, \fz]$.
It is easy to see that $H^{p,\,q}(\cx)\subset H_\cd^{p,\,q}(\cx)$.
\end{rem}

Using Remark \ref{r2.6}, we have the following conclusion.

\begin{lem}\label{l2.5}
Let $\cd$ be a collection of balls in $\cx$, $p\in(0,1]$ and $q\in
[1,\fz)\cap(p,\fz)$. Then $H_\cd^{p,\,q}(\cx)=H_\cd^{p,\,\fz}(\cx)$
with equivalent norms.
\end{lem}

\begin{pf}\rm
Notice that $(p,\,\fz)$-atoms and $(p,\,\fz)_\cd$-atoms are
$(p,\,q)$-atoms and $(p,\,q)_\cd$-atoms, respectively. Then from
Definition \ref{d2.7}, it follows that $H_\cd^{p,\,\fz}(\cx)\subset
H_\cd^{p,\,q}(\cx)$.

Conversely, let $f\in H_\cd^{p,\,q}(\cx)$. Then by Definition
\ref{d2.7}, there exist $\{\lz_j\}_{j\in\nn}$,
$\{\nu_k\}_{k\in\nn}\subset\cc$, $(p,\,q)$-atoms $\{a_j\}_{j\in\nn}$
and  $(p,\,q)_\cd$-atoms $\{b_k\}_{k\in\nn}$ such that
$$f=\sum_{j\in\nn}\lz_ja_j+\sum_{k\in\nn}\nu_kb_k,$$
which converges in $L^1(\cx)$ when $p=1$ or in
$(\llip(1/p-1;\,\cx))^\ast$ when $p\in(0,1]$, and
\begin{equation}\label{2.6}
\sum_{j\in\nn}|\lz_j|^p+\sum_{k\in\nn}|\nu_k|^p
\ls\|f\|^p_{H_\cd^{p,\,q}(\cx)}.
\end{equation}
For $k\in\nn$, assume that $\supp b_k\subset B_k\in\cd$ and let
$c_k\equiv[b_k-(b_k)_{B_k}\chi_{B_k}]/2$. Then it follows from
 Definition \ref{d2.5} that there exists a positive constant $\wz C$ such that
 $\{\wz Cc_k\}_{k\in\nn}$ are $(p,\,q)$-atoms,
$\{(b_k)_{B_k}\chi_{B_k}\}_{k\in\nn}$ are $(p,\,\fz)_\cd$-atoms
 and $b_k=2c_k+(b_k)_{B_k}\chi_{B_k}$.
 Moreover,
 $$f=\sum_{j\in\nn}\lz_ja_j+
 \sum_{k\in\nn}2\nu_kc_k+\sum_{k\in\nn}\nu_k(b_k)_{B_k}\chi_{B_k}, $$
 which converges in $L^1(\cx)$ when $p=1$ or in $(\llip(1/p-1;\,\cx))^\ast$
 when $p\in(0, 1)$.
By Remark \ref{r2.4} (ii) and \eqref{2.6}, we see that
$\sum_{j\in\nn}\lz_ja_j+2\sum_{k\in\nn}\nu_kc_k$
also converges in $L^1(\cx)$ when $p=1$
or in $(\lip(1/p-1;\,\cx))^\ast$ when $p\in(0, 1)$.
Let $g\equiv\sum_{j\in\nn}\lz_ja_j+2
\sum_{k\in\nn}\nu_kc_k$. Then Definition \ref{d2.6} together with Remark \ref{r2.5}
implies that $g\in H^{p,\,q}(\cx)=H^{p,\,\fz}(\cx)$.
Form this, Remark \ref{r2.6} and \eqref{2.6}, we deduce that
$g \in H^{p,\,\fz}(\cx)\subset H_\cd^{p,\,\fz}(\cx)$ and
$$\|g\|_{H_\cd^{p,\,\fz}(\cx)}\ls
\|g\|_{H^{p,\,\fz}(\cx)}\ls\|g\|_{H^{p,\,q}(\cx)}
\ls\|f\|_{H_\cd^{p,\,q}(\cx)},$$
which further implies that $f\in H_\cd^{p,\,\fz}(\cx)$ and
by \eqref{2.6},
\begin{eqnarray*}
\|f\|_{H_\cd^{p,\,\fz}(\cx)}&&\ls \|g\|_{H_\cd^{p,\,\fz}(\cx)}+
\lf\|\sum_{k\in\nn}\nu_k(b_k)_{B_k}\chi_{B_k}\r\|_{H_\cd^{p,\,\fz}(\cx)}\\
&&\ls\|f\|_{H_\cd^{p,\,q}(\cx)}+
\lf\{\sum_{k\in\nn}|\nu_k|^p\r\}^{1/p}
\ls\|f\|_{H_\cd^{p,\,q}(\cx)}.
\end{eqnarray*}
This finishes the proof of Lemma \ref{l2.5}.
\end{pf}

\begin{rem}\label{r2.7}\rm
(i) Let $\cd$ be a collection of balls in $\cx$, $p\in(0,1]$ and
$q\in [1,\fz]\cap(p,\fz]$. In what follows, based on Lemma \ref{l2.5},
we denote $H_\cd^{p,\,q}(\cx)$ simply by $H_\cd^p(\cx)$.

(ii)  Let $L^\fz_b(\cx)$ be the set of all functions of $L^\fz(\cx)$
with bounded support. Then from Definitions \ref{d2.6} and
\ref{d2.7}, it follows that
 $L^\fz_b(\cx)\cap H_\cd^p(\cx)$ is dense in $H_\cd^p(\cx)$
 and  $L^\fz_b(\cx)\cap H^p(\cx)$ is dense in $H^p(\cx)$.
\end{rem}

\begin{thm}\label{t2.1}
Let $\cd$ be a collection of balls in $\cx$ and $p\in(0,1]$.
Then $\ce^{1/p-1,\,1}_\cd(\cx)=(H_\cd^p(\cx))^\ast$.
\end{thm}

\begin{pf} \rm
We first prove
$\ce_\cd^{1/p-1,\,1}(\cx)\subset(H_\cd^{p,\,\fz}(\cx))^\ast$ for
$p\in(0,\,1]$. Let $f\in \ce_\cd^{1/p-1,\,1}(\cx)$. For all
$(p,\,\fz)$-atoms $a$ supported in $B\notin\cd$,
by Definition \ref{d2.5} (i), we have
$$\begin{array}{cl}
\lf|\dint_\cx f(x)a(x)\,d\mu(x)\r|
&=\lf|\dint_\cx [f(x)-f_B]a(x)\,d\mu(x)\r|\\
&\le\dfrac1{[\mu(B)]^{1/p}}\dint_B |f(x)-f_B|\,d\mu(x) \le
\|f\|_{\ce_\cd^{1/p-1,\,1}(\cx)}.
\end{array}$$
For all $(p,\,\fz)_\cd$-atoms $b$ supported in $B\in\cd$,
we also obtain
$$\begin{array}{cl}
\lf|\dint_\cx f(x)b(x)\,d\mu(x)\r| \le\dfrac1{[\mu(B)]^{1/p}}\dint_B
|f(x)|\,d\mu(x) \le \|f\|_{\ce_\cd^{1/p-1,\,1}(\cx)}.
\end{array}$$

Let $N\in\nn$ and $f_N\equiv\max\{\min\{f,N\},-N\}$. We claim
that $f_N\in\ce_\cd^{1/p-1,\,1}(\cx)$ and
\begin{equation}\label{2.7}
\|f_N\|_{\ce_\cd^{1/p-1,\,1}(\cx)}\le\frac94\|f\|_{\ce_\cd^{1/p-1,\,1}(\cx)}.
\end{equation}
In fact, if $B\in\cd$, then
$$\frac{1}{[\mu(B)]^{1/p}}\int_B|f_N(x)|\,d\mu(x)\le
\frac{1}{[\mu(B)]^{1/p}}\int_B|f(x)|\,d\mu(x)
\le\|f\|_{\ce_\cd^{1/p-1,\,1}(\cx)}.$$
Let $B\notin\cd$. For all $f$, $h\in \ce_\cd^{1/p-1,\,1}(\cx)$
and $g\equiv\max\{f,\,h\}$, we have that $g=(f+h+|f-h|)/2$ and
\begin{eqnarray*}
&&\frac{1}{[\mu(B)]^{1/p}}
\int_B|g(x)-g_B|\,d\mu(x)\\
&&\hs\le\frac{1}{2[\mu(B)]^{1/p}} \int_B|f(x)-f_B|\,d\mu(x)+
\frac{1}{2[\mu(B)]^{1/p}}
\int_B|h(x)-h_B|\,d\mu(x)\\
&&\hs\hs+\frac{1}{[\mu(B)]^{1/p}}\int_B|(f-h)(x)-(f-h)_B|\,d\mu(x)\\
&&\hs\le\frac32(\|f\|_{\ce^{1/p-1,\,1}(\cx)}+
\|h\|_{\ce^{1/p-1,\,1}(\cx)}).
\end{eqnarray*}
Similarly, for all $B\notin\cd$, $f$, $h\in \ce_\cd^{1/p-1,\,1}(\cx)$
and $\wz g\equiv\min\{f,\,h\}$, we have
\begin{eqnarray*}
&&\frac{1}{[\mu(B)]^{1/p}} \int_B|\wz g(x)-{\wz g}_B|\,d\mu(x)
\le\frac32\lf(\|f\|_{\ce^{1/p-1,\,1}(\cx)}+
\|h\|_{\ce^{1/p-1,\,1}(\cx)}\r).
\end{eqnarray*}
If $h\equiv N$ or $h\equiv-N$, then
$\|h\|_{\ce^{1/p-1,\,1}(\cx)}=0$.
By these facts and the definition of $f_N$, we have that
for all $B\notin \cd$,
$$\frac{1}{[\mu(B)]^{1/p}} \int_B|f_N(x)-(f_N)_B|\,d\mu(x)\le
\frac94\|f\|_{\ce_\cd^{1/p-1,\,1}(\cx)},$$
which implies the claim.

For all $g\in L^\fz_b(\cx)\cap H_\cd^{p,\,\fz}(\cx)$, since $fg\in
L^1(\cx)$, we define
$\ell(g)\equiv\int_\cx f(x)g(x)\,d\mu(x)$
and $\ell_N(g)\equiv\int_\cx f_N(x)g(x)\,d\mu(x).$
Moreover, there exist $\{\lz_j\}$, $\{\nu_k\}_{k\in\nn}\subset\cc$,
$(p,\,\fz)$-atoms $\{a_j\}_{j\in\nn}$  and  $(p,\,\fz)_\cd$-atoms
$\{b_k\}_{k\in\nn}$ such that
$$g=\sum_{j\in\nn}\lz_ja_j+\sum_{k\in\nn}\nu_kb_k$$
which converges in $L^1(\cx)$ when $p=1$ or in
$(\llip(1/p-1;\,\cx))^\ast$ when $p\in(0, 1)$, and
\begin{equation}\label{2.8}
\sum_{j\in\nn}|\lz_j|^p+\sum_{k\in\nn}|\nu_k|^p
\le2\|g\|^p_{H_\cd^{p,\,\fz}(\cx)}.
\end{equation}
By $f_N\in \ce_\cd^{1/p-1,\,1}(\cx)$ and $g\in H_\cd^{p,\,\fz}(\cx)$, we have
$$\ell_N(g)=\sum_{j\in\nn} \int_\cx f_N(x)\lz_ja_j(x)\,d\mu(x)
+\sum_{k\in\nn}\int_\cx f_N(x)\nu_kb_k(x)\,d\mu(x),$$
from which together with \eqref{2.7}, \eqref{2.8} and Remark \ref{r2.4} (ii),
 it follows that
$$|\ell_N(g)|\ls\|f_N\|_{\ce_\cd^{1/p-1,\,1}(\cx)}
\lf\{\dsum_{j\in\nn}|\lz_j|+\sum_{k\in\nn}|\nu_k|\r\}
\ls\|f\|_{\ce_\cd^{1/p-1,\,1}(\cx)}\|g\|_{H_\cd^{p,\,\fz}(\cx)}.$$
By this and the Lebesgue dominated theorem, we have
$$|\ell(g)| =\dlim_{N\rightarrow\fz}\lf|\dint_\cx
f_N(x)g(x)\,d\mu(x)\r|
\ls\|f\|_{\ce_\cd^{1/p-1,\,1}(\cx)}\|g\|_{H_\cd^{p,\,\fz}(\cx)},$$
which together with the density of $L^\fz_b(\cx)\cap
H_\cd^{p,\,\fz}(\cx)$ in $H_\cd^{p,\,\fz}(\cx)$ (see Remark
\ref{r2.7} (ii)) implies that $\ell\in (H_\cd^{p,\,\fz}(\cx))^\ast$
and
$\|\ell\|_{(H_\cd^{p,\,\fz}(\cx))^\ast}\ls\|f\|_{\ce_\cd^{1/p-1,\,1}(\cx)}.$
Thus,
\begin{equation}\label{2.9}
\ce_\cd^{1/p-1,\,1}(\cx)\subset (H_\cd^p(\cx))^\ast.
\end{equation}

We now prove that $(H_\cd^{p,\,2}(\cx))^\ast\subset
\ce_\cd^{1/p-1,\,2}(\cx)$. Let $\ell\in(H_\cd^{p,\,2}(\cx))^\ast$.
Since $H^{p,\,2}(\cx)\subset H_\cd^{p,\,2}(\cx)$, then
 $\ell\in (H^{p,\,2}(\cx))^\ast=\ce^{1/p-1,\,2}(\cx)$ (see Remark
 \ref{r2.5} and Remark \ref{2.1} (i)  and (iii)).
 Hence there exists
 $\wz f\in \ce^{1/p-1,\,2}(\cx)$ such that for
all constants $C$ and $g\in L^2(\cx)$ satisfying that
$\int_\cx g(x)\,d\mu(x)=0$ and $\supp(g)$ is bounded,
\begin{equation}\label{2.10}
\ell(g)=\int_\cx\wz f(x)g(x)\,d\mu(x)=\int_\cx(\wz f(x)+C)g(x)\,d\mu(x),
\end{equation}
and $\|\wz f\|_{\ce^{1/p-1,\,2}(\cx)}\ls\|\ell\|_{(H^{p,\,2}(\cx))^\ast}
\ls\|\ell\|_{(H^{p,\,2}_\cd(\cx))^\ast}$.
We then need to choose a suitable constant $C$ such that $f\equiv\wz f+C\in
\ce_\cd^{1/p-1,\,2}(\cx)$.

Observe that for all constants $\wz C$, $\wz f+\wz C\in  \ce^{1/p-1,\,2}(\cx)$.
Then by Lemma \ref{l2.1} (i), to show $f\in \ce_\cd^{1/p-1,\,2}(\cx)$ and
$\|f\|_{\ce_\cd^{1/p-1,\,2}(\cx)}\ls \|\ell\|_{(H^{p,\,2}_\cd(\cx))^\ast}$,
it suffices to show that for all $B\in \cd$,
\begin{equation}\label{2.11}
|f_B|[\mu(B)]^{1-1/p}\ls \|\ell\|_{\lf(H^{p,\,2}_\cd(\cx)\r)^\ast}.
\end{equation}
To this end, for any $B\in\cd$, let $L^2(B)\equiv\{f\in L^2(\cx):\,
\supp(f)\subset B\}$ and $L^2_0(B)
\equiv\{f\in L^2(B):\, \int_\cx f(x)\,d\mu(x)=0\}$.
Then for any $g\in L^2(B)$, the function
$g[\mu(B)]^{1/2-1/p}\|g\|^{-1}_{L^2(B)}$ is a $(p,\,2)_\cd$-atom
supported in $B$ and
$$|\ell(g)|\le\|\ell\|_{\lf(H^{p,\,2}_\cd(\cx)\r)^\ast}\|g\|_{H_\cd^{p,\,2}(\cx)}
\le[\mu(B)]^{1/p-1/2}\|\ell\|_{\lf(H^{p,\,2}_\cd(\cx)\r)^\ast}\|g\|_{L^2(B)},$$
which implies that
$\ell\in (L^2(B))^\ast=L^2(B)$. By this together
with the Riesz representation theorem,
there exists a function $f^B\in L^2(B)$
such that for all $g\in L^2(B)$, $\ell(g)=\int_B f^B(x)g(x)\,d\mu(x)$ and
\begin{equation}\label{2.12}
\lf\|f^B\r\|_{L^2(B)}\le [\mu(B)]^{1/p-1/2}\|\ell\|_{\lf(H^{p,\,2}_\cd(\cx)\r)^\ast}.
\end{equation}
Moreover, from this fact and \eqref{2.10},
we deduce that for all $g \in L^2_0(B)$,
$\int_\cx [f^B(x)-\wz f(x)]g(x)\,d\mu(x)=0$, which
further implies that $f^B-\wz f=0$ in $[L^2_0(B)]^\ast$.
Recall that $[L^2_0(B)]^\ast=L^2(B)/\cc$ (the space of functions
$f\in L^2(B)$ modulo constant functions) and
$f=0$ in $L^2(B)/\cc$ if and only if $f$
is a constant (see \cite[p.\,633]{cw77}).
Using these facts, we have that
$f^B-\wz f$ is a constant $C_B$.

Now it suffices to verify that for all balls $B,\,S\in\cd$, we have
$C_B=C_S$. Observe that
$g\equiv\{[\mu(\frac12 B)]^{-1}\chi_{\frac12 B}
-[\mu(\frac12 S)]^{-1}\chi_{\frac12 S}\}$ is a multiple of
certain $(p,\,2)$-atom, and $[\mu(\frac12 B)]^{-1}\chi_{\frac12 B}$ and
$\mu(\frac12 S)]^{-1}\chi_{\frac12 S}$ are multiplies of $(p,\,2)_\cd$-atoms.
Therefore, from the fact that
$f^B-C_B=\wz f=f^S-C_S$ and \eqref{2.10}, it follows that
\begin{eqnarray*}
\ell(g)&&=\ell\lf(\lf[\mu\lf(\frac12 B\r)\r]^{-1}\chi_{\frac12 B}\r)
-\ell\lf(\lf[\mu\lf(\frac12 S\r)\r]^{-1}\chi_{\frac12 S}\r)\\
&&=\dfrac1{\mu(\frac12 B)}\int_B f^B(x)\chi_{\frac12 B}(x)\,d\mu(x)-
\dfrac1{\mu(\frac12 S)}\int_S f^S(x)\chi_{\frac12 S}(x)\,d\mu(x)\\
&&=\int_{B\cup S}\wz f(x) g(x)\,d\mu(x)+C_B-C_S=\ell(g)+C_B-C_S,
\end{eqnarray*}
which implies that $C_B=C_S$. Denote the constant as above  by $\wz
C$ and define $f\equiv\wz f+\wz C$. Then by this, \eqref{2.12}
and the H\"older inequality,
we have that for all $B\in\cd$,
$$|f_B|[\mu(B)]^{1-1/p}=|(f^B)_B|[\mu(B)]^{1-1/p}
\ls\|\ell\|_{\lf(H^{p,\,2}_\cd(\cx)\r)^\ast}.$$
This implies \eqref{2.11}, from which and Lemma \ref{l2.1} (i),
we further deduce that $f\in \ce_\cd^{1/p-1,\,2}(\cx)$ and
$\|f\|_{\ce_\cd^{1/p-1,\,2}(\cx)}\ls\|\ell\|_{(H^{p,\,2}_\cd(\cx))^\ast}$.
Thus, $(H_\cd^p(\cx))^\ast\subset \ce_\cd^{1/p-1,\,2}(\cx)$,
which together with Lemma \ref{l2.2} and \eqref{2.9}
then completes the proof of Theorem \ref{t2.1}.
\end{pf}

\section{Boundedness of the radial and the Poisson maximal functions}\label{s3}

\hskip\parindent This section is devoted to the boundedness
of the radial and the Poisson maximal functions from $\ce^{\az,\,p}_{\rho}(\cx)$
to $\wz \ce^{\az,\,p}_{\rho}(\cx)$. We start with the notion
of the radial maximal function.

\begin{defn}\label{d3.1}\rm
Let $\rho$ be an admissible function
on $\cx$ and $\{T_t\}_{t>0}$ a family of linear integral
operators on $L^2(\cx)$. Moreover,
assume that there exist positive constants $C$, $\gz$, $\dz_1$, $\dz_2$,
$\bz$ satisfying that for all $t\in(0,\,\fz)$ and
$x,\,x',\,y\in\cx$ with $d(x,\,x')\le t/2$,
\begin{equation}\label{3.1}
|T_t(x,\,y)|\le C\frac1{V_t(x)+V(x,\,y)}\lf(\frac
t{t+d(x,\,y)}\r)^{\gz}\lf(\frac{\rho(x)}{t+\rho(x)}\r)^{\dz_1};
\end{equation}
\begin{equation}\label{3.2}
|T_t(x,\,y)-T_t(x',\,y)|\le C\frac1{V_t(x)+V(x,\,y)}\lf(\frac
t{t+d(x,\,y)}\r)^{\gz}\lf(\frac{d(x,\,x')}{t}\r)^{\bz};
\end{equation}
\begin{equation}\label{3.3}
|1-T_t(1)(x)|\le C\lf(\frac{t}{t+\rho(x)}\r)^{\dz_2}.
\end{equation}
\end{defn}

Let $\{T_t\}_{t>0}$ be as in Definition \ref{d3.1}.
For all $f\in L^1_\loc(\cx)$, the radial maximal function
$T^+$ is defined by
$$T^+(f)\equiv\sup_{t>0}|T_t(f)|.$$

Then we have the following result.

\begin{thm}\label{t3.1}
Let $\az\in(-\fz,\,\gz/n)\cap(-\fz,\,\min\{\bz/(2n),\,\dz_1/n,\,\dz_2/(2n)\} ]$,
$p\in(1,\,\fz)$ and $\rho$ be an admissible function.
If $\{T_t\}_{t>0}$ satisfies \eqref{3.1} through \eqref{3.3},
then there exists a positive
constant $C$ such that for all $f\in\ce^{\az,\,p}_{\rho}(\cx)$,
$T^+(f)\in\wz\ce^{\az,\,p}_\rho(\cx)$ and
$$\|T^+(f)\|_{\wz\ce^{\az,\,p}_{\rho}(\cx)}
\le C\|f\|_{\ce^{\az,\,p}_{\rho}(\cx)}.$$
\end{thm}

\begin{pf}\rm
We only consider the case that $\az\in(0,\,\gz/ n)\cap(0,\,
\min\{\bz/(2n),\,\dz_1/n,\,\dz_2/(2n)\}]$,
the proof for $\az\in(-\fz, 0]$ is similar
but easier. By the homogeneity of
$\|\cdot\|_{\ce^{\az,\,p}_{\rho}(\cx)}$ and $\|\cdot\|_{\wz
\ce^{\az,\,p}_{\rho}(\cx)}$, we assume that $f\in\ce^{\az,\,p}_{\rho}(\cx)$
and $\|f\|_{\ce^{\az,\,p}_{\rho}(\cx)}=1$.

Let $\cd_\rho$ be as in Remark \ref{r2.2} (iv) and
$B\equiv B(x_0,\,r)\in{\cd_\rho}$.
Observe that $T^+(f)\ls \hl(f)$,
where for all $x\in\cx$ and $f\in L^1_\loc(\cx)$,
$\hl(f)$ denotes the
Hardy-Littlewood maximal function of $f$ defined by
$$\hl(f)(x)\equiv\sup_{B\ni x}\frac1{\mu(B)}\dint_B|f(y)|\,d\mu(y).$$
Recall that $\hl$ is bounded on $\lp$ for $p\in(1, \fz]$.
Therefore $T^+$ is bounded on $\lp$ for all $p\in(1,\fz]$.
By this fact together with \eqref{2.1}, we see that
\begin{equation}\label{3.4}
\int_B [T^+(f\chi_{2B})(x)]^p\,d\mu(x) \ls
\int_{2B}|f(x)|^p\,d\mu(x)\ls[\mu(B)]^{1+\az p}.
\end{equation}
If $t\in(0, r)$, then by \eqref{3.1}, \eqref{2.1}, the H\"older inequality
and $\gz>\az n$, we have
\begin{eqnarray}\label{3.5}
\lf|T_t\lf(f\chi_{(2B)^\complement}\r)(x)\r|
&&\ls\int_{(2B)^\complement}\frac1{V_t(x)+V(x,\,y)}
\lf(\frac{t}{t+d(x,\,y)}\r)^{\gz}
|f(y)|\,d\mu(y)\\
&&\ls \sum_{j=1}^\fz 2^{-j\gz}\lf(\frac1{\mu(2^{j+1}B)}
\int_{2^{j+1}B}|f(y)|^p\,d\mu(y)\r)^{1/p}\nonumber\\
&&\ls \sum_{j=1}^\fz 2^{-j\gz} \lf[\mu\lf(2^{j+1}B\r)\r]^\az\ls
[\mu(B)]^\az\sum_{j=1}^\fz 2^{-j(\gz-\az n )}\ls [\mu(B)]^\az.\nonumber
\end{eqnarray}
Let $t\in[r, \fz)$. By \eqref{2.2}, we see that for all $a\in(0, \fz)$,
there exists a constant $\wz C_a\in[1, \fz)$ such that
for all $x$, $y\in\cx$ with $d(x, y)\le a\rho(x)$,
\begin{equation}\label{3.6}
\rho(y)/\wz {C}_a\le\rho(x)\le {\wz C}_a\rho(y).
\end{equation}
Recall that $B\in{\cd_\rho}$, which is equivalent
to that $r\ge \rho(x_0)$.
These facts imply that
for all $x\in B$, $\rho(x)\ls r$. By this together with
\eqref{3.1}, \eqref{2.1}, the H\"older inequality and
the facts that $\gz>\az n$ and $\dz_1\ge \az n$,
we have that for all $t\in[r,\fz)$ and $x\in B$,
\begin{eqnarray*}
\lf|T_t\lf(f\chi_{(2B)^\complement}\r)(x)\r|
&&\ls\int_{(2B)^\complement}\frac{|f(y)|}{V_t(x)+V(x,\,y)}
\lf(\frac{t}{t+d(x,\,y)}\r)^{\gz}\lf(\frac{\rho(x)}{t+\rho(x)}\r)^{\dz_1}\,d\mu(y)\\
&&\ls\lf(\frac{\rho(x)}{t+\rho(x)}\r)^{\dz_1}\sum_{j=1}^\fz
2^{-j\gz}\frac1{V_{2^{j-1}t}(x)}
\int_{d(x,\,y)<2^jt}|f(y)|\,d\mu(y)\nonumber\\
&&\ls\lf(\frac{\rho(x)}{t+\rho(x)}\r)^{\dz_1}\sum_{j=1}^\fz2^{-j\gz}\lf(\frac1{V_{2^{j+1}t}(x_0)}
\int_{d(x_0,\,y)<2^{j+1}t}|f(y)|^p\,d\mu(y)\r)^{1/p}\nonumber\\
&&\ls \lf(\frac{\rho(x)}{t+\rho(x)}\r)^{\dz_1}\sum_{j=1}^\fz2^{-j\gz}[V_{2^{j+1}t}(x_0)]^\az\nonumber\\
&&\ls\lf(\frac{\rho(x)}{t+\rho(x)}\r)^{\dz_1}\lf(\frac tr\r)^{\az n
}[\mu(B)]^\az\sum_{j=1}^\fz 2^{-j(\gz-\az n )} \ls [\mu(B)]^\az.\nonumber
\end{eqnarray*}
Combining this and \eqref{3.5} yields that for all $t\in(0, \fz)$,
$$\int_B \lf[T^+\lf(f\chi_{(2B)^\complement}\r)(x)\r]^p\,d\mu(x)
\ls [\mu(B)]^{1+\az p},$$
which together with \eqref{3.4} gives us that
$$\int_B [T^+(f)(x)]^p\,d\mu(x)
\ls [\mu(B)]^{1+\az p}.$$
This also implies that $T^+(f)(x)<\fz$ for
$\mu$-a.\,e. $x\in\cx$.

It remains to show that for all $B\equiv B(x_0,\,r)\notin{\cd_\rho}$,
\begin{equation*}
\int_B\lf[T^+(f)(x)-{\mathop\einf_B}T^+(f)
\r]^p\,d\mu(x)\ls[\mu(B)]^{1+\az p}.
\end{equation*}
Let $f_1\equiv(f-f_B)\chi_{2B}$,
$f_2\equiv(f-f_B)\chi_{(2B)^\complement}$, $B_1\equiv\lf\{x\in B:\
T_{r}^+(f)(x)\ge T_\fz^+(f)(x)\r\}$ and $B_2\equiv B\setminus B_1$,
where $T_{r}^+(f)\equiv\sup_{0<t<4r}|T_t(f)|$ and
$T_{\fz}^+(f)\equiv\sup_{t\ge4r}|T_t(f)|.$ We have
\begin{eqnarray*}
&&\int_B\lf[T^+(f)(x)-{\mathop\einf_{B}}T^+(f)\r]^p\,d\mu(x)\\
&&\hs\ls
\int_{B_1}\lf[T_{r}^+(f)(x)-{\mathop\einf_B}|T_r(f)|\r]^p\,d\mu(x)\\
&&\hs\hs +\int_{B_2}\lf[T_{\fz}^+(f)(x)-
{\mathop\einf_B}T_{\fz}^+(f)\r]^p\,d\mu(x)\\
&&\hs\ls \int_B[T_r^+(f_1)(x)]^p\,d\mu(x)
+\mu(B)\sup_{x,\,y\in B}\sup_{0<t<4r}\lf|T_t(f_B)(x)-T_r(f)(y)\r|^p\\
&&\hs\hs+\int_B\lf[T_r^+(f_2)(x)\r]^p \,d\mu(x)
+\mu(B)\sup_{x,\,y\in B}\sup_{t\ge4r}|T_t(f)(x)-T_t(f)(y)|^p\\
&&\hs\equiv{\rm E_1}+{\rm E_2}+{\rm E_3}+{\rm E_4}.
\end{eqnarray*}
By the H\"older inequality, $L^p(\cx)$-boundedness
of $T^+$ and \eqref{2.1}, we have
$${\rm E_1}\ls\int_{2B}|f(x)-f_B|^p\,d\mu(x)\ls[\mu(B)]^{1+\az p}.$$
On the other hand, using \eqref{3.1}, \eqref{2.1},
the H\"older inequality, Lemma \ref{l2.4} (ii)
and $\gz>\az n$, we have
that for all $t\in(0, 4r)$ and $x\in B$,
\begin{eqnarray*}
|T_t(f_2)(x)| &&\ls\int_{(2B)^\complement}\frac1{V_t(x)+V(x,\,z)}
\lf(\frac{t}{t+d(x,\,z)}\r)^{\gz}
|f(z)-f_B|\,d\mu(z)\\
&&\ls\sum_{j=1}^\fz 2^{-j\gz}
\frac1{V_{2^{j-1}r}(x)}\int_{2^{j+1}B}[|f(z)-
f_{2^{j+1}B}|+|f_B-f_{2^{j+1}B}|]\,d\mu(z)\\
&&\ls[\mu(B)]^\az\sum_{j=1}^\fz 2^{-j(\gz-\az n )}\ls[\mu(B)]^\az.
\end{eqnarray*}
This implies that ${\rm E_3}\ls[\mu(B)]^{1+\az p}$.

Similarly, by applying \eqref{3.1},
\eqref{2.1} and $\gz>\az n$,
we have that for all $x\in B$,
\begin{eqnarray}\label{3.7}
|T_r(f-f_B)(x)|&&\ls\int_\cx\frac1{V_r(x)+V(x,\,z)}\lf(\frac
r{r+d(x,\,z)}\r)^{\gz}
|f(z)-f_B|\,d\mu(z)\\
&&\ls\sum_{j=0}^\fz 2^{-j\gz}
\frac1{V_{2^{j-1}r}(x)}\int_{2^{j+1}B}|f(z)-f_B|\,d\mu(z)\ls[\mu(B)]^\az.\nonumber
\end{eqnarray}
From Lemma \ref{l2.4} (i), \eqref{3.3}, $\dz_2\ge\az n$ and $t<4r\ls \rho(x_0)$,
it follows that for all $x\in B$,
\begin{eqnarray*}
|f_B-T_t(f_B)(x)|&&=|f_B| |1-T_t(1)(x)| \ls[\mu(B)]^\az\lf(\frac
t{\rho(x_0)}\r)^{\dz_2-\az n } \ls[\mu(B)]^\az.
\end{eqnarray*}
This together with \eqref{3.7} implies that
\begin{eqnarray*}
{\rm E_2}&\ls&\mu(B)\sup_{x,\,y\in B}\sup_{0<t<4r}\lf\{|T_t(f_B)(x)-f_B|^p+
|f_B-T_r(f_B)(y)|^p+|T_r(f_B-f)(y)|^p\r\}\\
&\ls&[\mu(B)]^{1+\az p}.
\end{eqnarray*}

To estimate ${\rm E_4}$, we first observe that for all $x,\,y\in B$,
$\rho(x)\sim\rho(x_0)\sim\rho(y)$ (see \eqref{3.6}).
By this and \eqref{3.2}, we have that for all $x$, $y\in B$ and $t\in [4r, \fz)$,
$$|T_t(1)(x)-T_t(1)(y)|\ls \lf(\frac{r}{t}\r)^\bz.$$
On the other hand, it follows from Lemma \ref{l2.4} (i) and \eqref{2.1} that
$$|f_{B(x_0,\,t)}|\ls \lf(\frac{\rho(x_0)}{r}\r)^{\az n}[\mu(B)]^\az.$$
Then by these facts and $\az n\le\min\{\frac{\bz}2, \frac{\dz_2}2\}$,
we obtain that for all $t\in[4r, \fz)$,
\begin{eqnarray*}
&&|T_t(1)(x)-T_t(1)(y)||f_{B(x_0,\,t)}|\\
&&\hs\ls\lf(\frac{\rho(x_0)}{r}\r)^{\az n}[\mu(B)]^\az
|T_t(1)(x)-T_t(1)(y)|^{\frac12}
\lf[|T_t(1)(x)-1| +|1-T_t(1)(y)|\r]^{\frac12}\\
&&\hs\ls\lf(\frac{\rho(x_0)}{r}\r)^{\az n}[\mu(B)]^\az
\lf(\frac{r}{\rho(x_0)}\r)^{\min\lf\{\frac{\bz}2,
\frac{\dz_2}2\r\}}\ls[\mu(B)]^\az.
\end{eqnarray*}
On the other hand, by \eqref{3.2}, \eqref{2.1},
the H\"older inequality, Lemma \ref{l2.4} (ii), $\gz>\az n$
and $\bz\ge \az n$, we see that
for all $x$, $y\in B$ and $t\in[4r, \fz)$,
\begin{eqnarray*}
&&\lf|T_t\lf(f-f_{B(x_0,\,t)}\r)(x)-T_t\lf(f-f_{B(x_0,\,t)}\r)(y)\r|\\
&&\hs\ls\int_\cx\lf(\frac {d(x,\,y)}{t}
\r)^{\bz}\frac1{V_t(x)+V(x,\,z)}
\lf(\frac t{t+d(x,\,z)}\r)^{\gz}|f(z)-f_{B(x_0,\,t)}|\,d\mu(z)\\
&&\hs\ls\lf(\frac rt\r)^{\bz}\sum_{j=0}^\fz
\frac{2^{-j\gz}}{V_{2^{j-1}t}(x)}\dint_{d(x,\,z)<2^jt}
\lf[|f(z)-f_{B(x_0,\,2^{j+1}t)}|+
\lf|f_{B(x_0,\,t)}-f_{B(x_0,\,2^{j+1}t)}\r|\r]\,d\mu(z)\\
&&\hs\ls\lf(\frac rt\r)^{\bz}\sum_{j=0}^\fz 2^{-j(\gz-\az n)}[\mu(B(x_0,\,t))]^\az
\ls[\mu(B)]^\az.
\end{eqnarray*}
These inequalities above lead to that
\begin{eqnarray*}
{\rm E_4}&&\ls\mu(B)\sup_{x,\,y\in B}\sup_{t\ge 4r}
\lf|T_t(f-f_{B(x_0,\,t)})(x)-T_t(f-f_{B(x_0,\,t)})(y)\r|^p\\
&&\hs+\mu(B)\sup_{x,\,y\in B}\sup_{t\ge 4r}
\lf[|T_t(1)(x)-T_t(1)(y)|\lf|f_{B(x_0,\,t)}\r|\r]^p
\ls[\mu(B)]^{1+\az p},
\end{eqnarray*}
which completes the proof of Theorem \ref{t3.1}.
\end{pf}

Now we consider the boundedness of the Poisson semigroup maximal operator.
Let $\{T_t\}_{t>0}$ be a family of linear
integral operators on $L^2(\cx)$. We always set
$$P_t\equiv\frac1{\sqrt\pi}\int_0^\fz \frac{e^{-s}}{\sqrt
s}T_{t/(2\sqrt s)}\,ds.$$ For all $f\in L_\loc^1(\cx)$, define
the Poisson semigroup maximal operator $P^+$ by
$$P^+(f)\equiv\sup_{t>0}|P_t(f)|.$$

\begin{lem}\label{l3.1}
Assume that $\{T_t\}_{t>0}$ satisfies \eqref{3.1} through  \eqref{3.3}
with the same constants $\dz_1,\,\dz_2,\,\bz,\,\gz$ as there.
Then $\{P_t\}_{t>0}$
also satisfies
 \eqref{3.1} through \eqref{3.3}
 with the constants $\dz_1$, $\dz'_2$, $\bz'$ and $\gz'$,  where
 $\dz'_2\in(0,\,1)\cap(0,\,\dz_2]$, $\bz'\in(0,\,1)\cap(0,\,\bz]$
 and $\gz'\in(0,\,1)\cap(0,\,\gz]$.
 \end{lem}

\begin{pf}\rm For all $a$, $s$, $t\in(0,\fz)$, from the fact that
$t+a\le (1+s)(t/s+a),$
it follows that
\begin{equation}\label{3.8}
\frac{t/s}{t/s+a}\le (1+s^{-1})\frac {t}{t+a}.
\end{equation}
On the other hand, from \eqref{2.1}, we deduce that
for all $x$, $y\in \cx$ and $s$, $t\in(0,\fz)$,
\begin{equation}\label{3.9}
\begin{array}[t]{ccl}
V_{t/s}(x)+V(x,\,y)&&\sim \mu(B(x,\,t/s+d(x,\,y)))\\
&&\gs(1+s)^{- n } \mu(B(x,\,t+d(x,\,y)))\sim (1+s)^{- n
}[V_{t}(x)+V(x,\,y)].
\end{array}
\end{equation}
By \eqref{3.1}, \eqref{3.8} and \eqref{3.9}, we see that for all $x$, $y\in \cx$,
\begin{eqnarray*}
\lf| P_t(x,\,y)\r|&&\ls\int_0^\fz e^{-s^2/4}  T_{t/s}(x,\,y)\,ds\\
&&\ls \int_0^\fz  e^{-s^2/4} \frac1{V_{t/s}(x)+V(x,\,y)}\lf(\frac
{t/s}{t/s+d(x,\,y)}\r)^{\gz} \lf(\frac{\rho(x)}{t/s+\rho(x)}\r)^{\dz_1}\,ds\nonumber\\
&&\ls \frac1{V_{t}(x)+V(x,\,y)}\lf(\frac{t}{t+d(x,\,y)}\r)^{\gz'}
\lf(\frac{\rho(x)}{t+\rho(x)}\r)^{\dz_1}\nonumber\\
&&\hs\times\int_0^\fz  e^{-s^2/4}(1+s)^{n+\dz_1} (1+s^{-\gz'})\,ds\nonumber\\
&&\ls\frac1{V_{t}(x)+V(x,\,y)}\lf(\frac{t}{t+d(x,\,y)}\r)^{\gz'}
\lf(\frac{\rho(x)}{t+\rho(x)}\r)^{\dz_1}.\nonumber
\end{eqnarray*}

Now we prove that for all $t\in(0, \fz)$ and $x$, $x'$, $y\in\cx$ with
$d(x,\,x')\le t/2$,
\begin{equation}\label{3.10}
\lf| P_t(x,\,y)- P_t(x',\,y)\r| \ls \lf(\frac
{d(x,\,x')}{t}\r)^{\bz'}\frac1{V_{t}(x)+V(x,\,y)}\lf(\frac
{t}{t+d(x,\,y)}\r)^{\gz'}.
\end{equation}
Observe that in this case,
$t+d(x,\,y)\sim t+d(x',\,y)$ and
$d(x,\,x')\le t/(2s)$ if and only if $s\le t/[2d(x,\,x')]$.  Then
\eqref{3.1} and \eqref{3.2} together with
\eqref{3.8} and \eqref{3.9} yield that
\begin{eqnarray*}
&&\lf| P_t(x,\,y)- P_t(x',\,y)\r|\\
&&\hs\ls\int_0^\fz e^{-s^2/4} \lf|
 T_{t/s}(x,\,y)- T_{t/s}(x',\,y)\r|\,ds\\
&&\hs\ls \lf[\int_0^{t/[2d(x,\,x')]}\lf(\frac
{d(x,\,x')}{t/s}\r)^{\bz}+\int_{t/[2d(x,\,x')]}^\fz\r]
\frac{e^{-s^2/4}}{V_{t/s}(x)+V(x,\,y)}\lf(\frac
{t/s}{t/s+d(x,\,y)}\r)^{\gz} \,ds\\
&&\hs\ls \lf[\int_0^{t/[2d(x,\,x')]} (1+s)^{\bz'}
+\int_{t/[2d(x,\,x')]}^\fz s^{\bz'}\r]e^{-s^2/4}(1+s)^n
(1+s^{-\gz'})\,ds\\
&&\hs\hs\times\lf(\frac
{d(x,\,x')}{t}\r)^{\bz'}\frac1{V_{t}(x)+V(x,\,y)}\lf(\frac
{t}{t+d(x,\,y)}\r)^{\gz'} \\
&&\hs\ls \lf(\frac
{d(x,\,x')}{t}\r)^{\bz'}\frac1{V_{t}(x)+V(x,\,y)}\lf(\frac
{t}{t+d(x,\,y)}\r)^{\gz'},
\end{eqnarray*}
which implies \eqref{3.10}.

On the other hand, by \eqref{3.3} and \eqref{3.8},
we see that for all $x\in\cx$ and $t\in(0, \fz)$,
\begin{eqnarray*}
\lf|1-P_t(1)(x)\r|&\ls& \int_0^\fz e^{-s^2/4}
\lf|1-T_{t/s}(1)(x)\r|\,ds\\
&\ls&\int_0^\fz e^{-s^2/4}\lf(\frac{t/s}{t/s+\rho(x)}\r)^{\dz_2}\,ds\\
&\ls&\lf(\frac{t}{t+\rho(x)}\r)^{\dz_2'} \int_0^\fz
e^{-s^2/4}(1+s^{-\dz_2'})\,ds
\ls\lf(\frac{t}{t+\rho(x)}\r)^{\dz_2'}.
\end{eqnarray*}
This finishes the proof of Lemma \ref{l3.1}.
\end{pf}

\begin{thm}\label{t3.2} Let $\rho$ be an admissible function
and $\{T_t\}_{t>0}$ satisfy \eqref{3.1}
through \eqref{3.3} with the same constants $\bz,\,\gz,\,\dz_1,\,\dz_2$
as there and $\dz'_2$, $\bz'$ and $\gz'$ be positive
constants such that $\dz'_2\in(0,\,1)\cap(0,\,\dz_2]$,
$\bz'\in(0,\,1)\cap(0,\,\bz]$
 and $\gz'\in(0,\,1)\cap(0,\,\gz]$.
Let $\az\in(-\fz,\,\gz'/n)\cap(-\fz,\,\min\{\bz'/(2n),
\,\dz_1/n,\,\dz'_2/(2n)\}]$ and $p\in(1,\,\fz)$.
Then there exists a positive constant $C$ such that
for all $f\in\ce^{\az,\,p}_{\rho}(\cx)$,
$P^+(f)\in\wz\ce^{\az,\,p}_\rho(\cx)$ and
$$\|P^+(f)\|_{\wz\ce^{\az,\,p}_{\rho}(\cx)}
\le C\|f\|_{\ce^{\az,\,p}_{\rho}(\cx)}.$$
\end{thm}

\begin{pf}\rm
Notice that our assumption on $\{T_t\}_{t>0}$ and Lemma \ref{l3.1}
imply that $\{P_t\}_{t>0}$ satisfies \eqref{3.1} through
\eqref{3.3} with constants $\dz_1$, $\dz'_2$, $\gz'$ and $\bz'$. By this and
an argument similar to the proof of Theorem \ref{t3.1}, we can prove
Theorem \ref{t3.2}. We omit the details by the similarity.
This finishes the proof of Theorem \ref{t3.2}.
\end{pf}

\begin{rem}\label{r3.1}\rm

(i) If $\az>0$, then by Lemma \ref{l2.2} (ii), the spaces
$\wz\ce^{\az,\,p}_\rho(\cx)$ in Theorems \ref{t3.1} and \ref{t3.2}
are exactly the spaces $\ce^{\az,\,p}_\rho(\cx)$.
 If $\az<0$ and $\cx$ is an RD-space, then by Remark \ref{r2.3} (ii)
and the fact that the maximal operators are nonnegative, we know that if the space
$\wz\ce^{\az,\,p}_\rho(\cx)$ in Theorems \ref{t3.1} and \ref{t3.2}
is replaced by the space $\ce^{\az,\,p}_\rho(\cx)$,
we obtain the same results.

(ii) Let $\cx$ be an RD-space and $\rho$ an admissible function.
Assume that there exist constants $C\in(0,\,\fz)$,  $\ez_1\in(0,1]$,
$\ez_2\in(0,\,\fz)$, $\dz\in(0,\,1]$ and $\gz\in(0,\,\fz)$, and an
$(\ez_1, \ez_2)$-$\ati$ $\{\wz T_t\}_{t>0}$ (see, for example,
\cite{hmy2, yyz} for the definition of $\ati$) with kernels $\{\wz
T_t(x, y)\}_{t>0}$ such that for all $t\in(0, \fz)$ and $x,\,
y\in\cx$,
\begin{equation}\label{3.11}
\lf|T_t(x,\,y)-\wz T_t(x,\,y)\r|\le
C\lf(\frac{t}{t+\rho(x)}\r)^{\dz} \frac1{V_t(x)+V(x,\,y)}\lf(\frac
t{t+d(x,\,y)}\r)^{\gz}.
\end{equation}
If $\az=0$ and \eqref{3.1} through \eqref{3.3}
were replaced by \eqref{3.11},
Theorems \ref{t3.1} and \ref{t3.2} were
obtained in \cite{yyz}. We remark that
since for all $x\in\cx$, $\wz T_t(1)(x)=1$ (see \cite{yyz}),
\eqref{3.11} implies \eqref{3.3} with
$\dz_2=\dz$.
\end{rem}

\section{Boundedness of the Littlewood-Paley
$g$-function \label{s4} }
\hskip\parindent In this section, we consider the boundedness
of certain variant of the Littlewood-Paley $g$-function from
$\ce^{\az,\,p}_{\rho}(\cx)$ to $\wz\ce^{\az,\,p}_{\rho}(\cx)$.
The boundedness from $\bmo_\rho(\cx)$ to $\blo_\rho(\cx)$ where
$\cx$ is an RD-space of this operator was obtained in \cite{yyz}.

Let $\rho$ be an admissible function on $\cx$ and $\{Q_t\}_{t>0}$
a family of  operators bounded on $\lt$ with  integral kernels
$\{Q_t(x,\,y)\}_{t>0}$ satisfying that there exist constants
$C\in(0,\,\fz)$, $\dz_1\in(0,\,\fz)$, $\dz_2\in(0,\,1)$, $\bz\in(0,\,1]$
 and $\gz\in(0,\,\fz)$ such that for all
$t\in(0,\,\fz)$ and $x,\,x',\,y\in\cx$ with $d(x,\,x')\le \frac t2$,

$(Q)_{\rm i}$ $|Q_t(x,\,y)|\le C\frac1{V_t(x)+V(x,\,y)}(\frac
t{t+d(x,\,y)})^{\gz}(\frac {\rho(x)}{t+\rho(x)})^{\dz_1}$;

$(Q)_{\rm ii}$  $|Q_t(x,\,y)-Q_t(x',\,y)|\le C(\frac
{d(x,\,x')}{t+d(x,\,y)})^\bz \frac{1}{V_t(x)+V(x,\,y)}(\frac
t{t+d(x,\,y)})^{\gz}$;

$(Q)_{\rm iii}$ $|\int_\cx Q_t(x,\,y)d\mu(y)|\le C(\frac
t{t+\rho(x)})^{\dz_2}$.

For all $f\in L^1_\loc(\cx)$ and $x\in\cx$, define the Littlewood-Paley
$g$-function by
\begin{equation}\label{4.1}
g(f)(x)\equiv\lf(\int_0^\fz |Q_t(f)(x)|^2\frac {dt}t\r)^{1/2}.
\end{equation}

\begin{lem}\label{l4.1}
Let $\az\in(-\fz,\,\min\{\gz/n,\,\dz_2/n\})$, $p\in(1,\fz)$ and
$\rho$ be an admissible function on $\cx$. Then there exists a
positive constant $C$ such that for all $f\in
\ce^{\az,\,p}_{\rho}(\cx)$,

(i) for all $x\in\cx$ and $t>0$,
$$|Q_t(f)(x)|\le C\lf(\frac{\rho(x)}{t+\rho(x)}\r)^{\dz_1}
[\mu(B(x,\,t))]^\az\|f\|_{\ce^{\az,\,p}_{\rho}(\cx)};$$

(ii) for all $x,\,y\in\cx$ and $t\ge2d(x,\,y)$,

$$|Q_t(f)(x)-Q_t(f)(y)|\le
\lf\{\begin{array}{ll}
C\lf(\frac{d(x,\,y)}{t}\r)^{\bz}\lf(1+\frac{\rho(x)}{t}\r)^{\az n }
[\mu(B(x,\,t))]^\az\|f\|_{\ce^{\az,\,p}_{\rho}(\cx)},&\hs \az>0;\\
C\lf(\frac{d(x,\,y)}{t}\r)^{\bz}\lf(1+\log\frac{\rho(x)}{t}\r)
[\mu(B(x,\,t))]^\az\|f\|_{\ce^{\az,\,p}_{\rho}(\cx)},&\hs\az\le 0.
\end{array}\r.$$
\end{lem}

\begin{pf}\rm
By the homogeneity of $\|\cdot\|_{\ce^{\az,\,p}_{\rho}(\cx)}$, we
may assume that $f\in \ce^{\az,\,p}_{\rho}(\cx)$ and
$\|f\|_{\ce^{\az,\,p}_{\rho}(\cx)}=1$.
By ${\rm (Q)_i}$, \eqref{4.2}, \eqref{2.1}, $\gz>\az n$ and the H\"older inequality,
we have that for all $x\in \cx$ and $t\ge\rho(x)$,
\begin{eqnarray}\label{4.2}
|Q_t(f)(x)|&&\ls\int_{\cx}\frac1{V_t(x)+V(x,\,y)}
\lf(\frac{t}{t+d(x,\,y)}\r)^{\gz}\lf(\frac{\rho(x)}{t+\rho(x)}\r)^{\dz_1}
|f(y)|\,d\mu(y)\\
&&\ls\lf(\frac{\rho(x)}{t+\rho(x)}\r)^{\dz_1}\sum_{j=0}^\fz
2^{-j\gz}\frac1{V_{2^{j-1}t}(x)}
\int_{d(x,\,y)<2^jt}|f(y)|\,d\mu(y)\nonumber\\
&&\ls \lf(\frac{\rho(x)}{t+\rho(x)}\r)^{\dz_1}\sum_{j=0}^\fz
2^{-j\gz} [\mu(B(x,\,2^jt))]^\az\nonumber\\
&&\ls\lf(\frac{\rho(x)}{t+\rho(x)}\r)^{\dz_1}[\mu(B(x,\,t))]^\az
\sum_{j=0}^\fz\max\lf\{2^{-j(\gz-\az n)},\ 2^{-j\gz}\r\}\nonumber\\
&&\ls\lf(\frac{\rho(x)}{t+\rho(x)}\r)^{\dz_1}[\mu(B(x,\,t))]^\az.\nonumber
\end{eqnarray}

Let $x\in\cx$ and $t<\rho(x)$. In this case, $t+\rho(x)\sim \rho(x)$.
Using $\gz>\az n$, ${\rm (Q)_i}$, \eqref{2.1},
Lemma \ref{l2.4} (ii) and the H\"older inequality, we have
\begin{eqnarray*}
&&\lf|Q_t\lf(f-f_{B(x,\,t)}\r)(x)\r|\\
&&\hs\ls\sum_{j=0}^\fz
2^{-j\gz}\frac1{V_{2^{j-1}t}(x)}
\int_{d(x,\,y)<2^jt}|f(y)-f_{B(x,\,t)}|\,d\mu(y)\\
&&\hs\ls\sum_{j=0}^\fz
2^{-j\gz}\lf\{\frac1{V_{2^{j}t}(x)}
\int_{d(x,\,y)<2^jt}\lf|f(y)-f_{B(x,\,2^jt)}\r|\,d\mu(y)
+\lf|f_{B(x,\,2^jt)}-f_{B(x,\,t)}\r|\r\}\\
&&\hs\ls [\mu(B(x,\,t))]^\az\sum_{j=0}^\fz 2^{-j\gz}
\max\lf\{2^{j\max\{\az n,\,0\}}, j+1\r\}\ls[\mu(B(x,\,t))]^\az.
\end{eqnarray*}

On the other hand, from ${\rm (Q)_{iii}}$, Lemma \ref{l2.4} (i), $t<\rho(x)$,
and the fact $\dz_2>\az n$, we deduce that
\begin{eqnarray*}
|Q_t(f_{B(x,\,t)})(x)|&&\ls[\mu(B(x,\,t))]^\az\lf(\frac
t{t+\rho(x)}\r)^{\dz_2} \max\lf\{1+\log\frac{\rho(x)}{t},
\lf(\frac {\rho(x)}t\r)^{\max\{\az n ,\,0\}}\r\}\\
&&\ls[\mu(B(x,\,t))]^\az\lf(\frac {\rho(x)}{t+\rho(x)}\r)^{\dz_1}.
\end{eqnarray*}
This gives (i).

To show (ii), by ${\rm (Q)_{ii}}$, we see that for all $x$, $y\in\cx$ and
$t\ge2d(x,\,y)$,
\begin{eqnarray}\label{4.3}
&&|Q_t(f)(x)-Q_t(f)(y)|\\
&&\hs\ls\int_{\cx}\lf(\frac{d(x,\,y)}
{t+d(x,\,z)}\r)^{\bz}\frac1{V_t(x)+V(x,\,z)}
\lf(\frac{t}{t+d(x,\,z)}\r)^{\gz}
|f(z)|\,d\mu(z)\nonumber\\
&&\hs\ls\lf(\frac{d(x,\,y)}{t}\r)^{\bz}\sum_{j=0}^\fz
2^{-j\gz}\frac1{V_{2^{j-1}t}(x)}
\int_{d(x,\,z)<2^jt}|f(z)|\,d\mu(z).\nonumber
\end{eqnarray}

Now we consider the following two cases.
{\it Case (i)} $\az\in(0,\fz)$. In this case, if $t\ge\rho(x)$, by $\gz> \az n$, the
H\"older inequality, \eqref{4.3} and \eqref{2.1}, we have
\begin{eqnarray}\label{4.4}
|Q_t(f)(x)-Q_t(f)(y)|&\ls&\lf(\frac{d(x,\,y)}{t}\r)^{\bz}\sum_{j=0}^\fz
2^{-j\gz}[\mu(B(x,\,2^jt))]^\az\\
&\ls&\lf(\frac{d(x,\,y)}{t}\r)^{\bz}[\mu(B(x,\,t))]^\az.\nonumber
\end{eqnarray}
Assume that $t<\rho(x)$. Let $N_1\in\nn$ such that
$2^{N_1-1}t<\rho(x)\le 2^{N_1}t$. From the H\"older inequality and
\eqref{2.1}, it follows that
\begin{eqnarray}\label{4.5}
&&\sum_{j=N_1}^\fz 2^{-j\gz}\frac1{V_{2^{j-1}t}(x)}
\int_{d(x,\,z)<2^jt}|f(z)|\,d\mu(z)\\
&&\hs\ls\sum_{j=N_1}^\fz2^{-j\gz}[\mu(B(x,\,2^jt))]^\az
\ls[\mu(B(x,\,t))]^\az.\nonumber
\end{eqnarray}
By the H\"older inequality, \eqref{2.1} and Lemma \ref{l2.4} (i), we see that
for all $j\in\{0, 1, \cdots, N_1-1\}$,
\begin{eqnarray*}
\frac1{V_{2^{j-1}t}(x)}
\int_{d(x,\,z)<2^jt}|f(z)|\,d\mu(z)
\ls \lf(\frac{\rho(x)}{t}\r)^{\az n}[\mu(B(x,\,t))]^{\az}.
\end{eqnarray*}
This together with $\gz> \az n$ gives us that
\begin{eqnarray*}
\sum_{j=0}^{N_1-1}2^{-j\gz}\frac1{V_{2^{j-1}t}(x)}
\int_{d(x,\,z)<2^jt}|f(z)|\,d\mu(z) \ls
\lf(\frac{\rho(x)}{t}\r)^{\az n}[\mu(B(x,\,t))]^{\az}.
\end{eqnarray*}
Combining this, \eqref{4.3} through \eqref{4.5} leads to that
for all $x$, $y\in\cx$ and $t\ge2d(x, y)$,
\begin{eqnarray*}
&&|Q_t(f)(x)-Q_t(f)(y)|\ls\lf(1+\frac{\rho(x)}{t}\r)^{\az
n}[\mu(B(x,\,t))]^{\az}.
\end{eqnarray*}

{\it Case (ii)} $\az\in(-\fz, 0]$. If $t\ge \rho(x),$ then \eqref{4.3} yields that
$$|Q_t(f)(x)-Q_t(f)(y)|\ls\lf(\frac{d(x,\,y)}{t}\r)^{\bz}[\mu(B(x,\,t))]^\az.$$
Let $t<\rho(x)$ and $N_1$ be the integer as in Case (i).
Then by \eqref{4.3}, \eqref{2.1}, Lemma \ref{l2.4} (i) and the H\"older inequality,
we have
\begin{eqnarray*}
&&|Q_t(f)(x)-Q_t(f)(y)|\\
&&\hs\ls\lf(\frac{d(x,\,y)}{t}\r)^{\bz}\lf\{\sum_{j=0}^{N_1-1}
2^{-j\gz}\frac1{V_{2^{j-1}t}(x)}
\int_{d(x,\,z)<2^jt}|f(z)|\,d\mu(z)+\sum_{j=N_1}^\fz\cdots\r\}\\
&&\hs\ls\lf(\frac{d(x,\,y)}{t}\r)^{\bz}\lf\{\sum_{j=0}^{N_1-1}
2^{-j\gz}\lf(1+\log\frac{\rho(x)}{t}\r)+\sum_{j=N_1}^\fz2^{-j\gz}\r\}[\mu(B(x,\,t))]^\az\\
&&\hs\ls\lf(\frac{d(x,\,y)}{t}\r)^{\bz}
\lf(1+\log\frac{\rho(x)}{t}\r)[\mu(B(x,\,t))]^\az,
\end{eqnarray*}
which implies (ii) and then completes the proof of Lemma \ref{l4.1}.
\end{pf}

\begin{thm}\label{t4.1}
Let $p\in(1,\,\fz)$, $\rho$ be an admissible function on $\cx$,
$g$ as in \eqref{4.1} and
$$\az\in(-\fz,\,\bz/(3 n )]\cap(-\fz,\,\min\{\gz/n,\,
\dz_1/ n ,\dz_2/(3 n) \}).$$
If $g(\cdot)$ is bounded on
$L^p(\cx)$, then there exists a positive constant $C$ such that for
all $f\in\ce^{\az,\,p}_{\rho}(\cx)$,
$[g(f)]^2\in\wz\ce^{2\az,\,p/2}_{\rho}(\cx)$ and
$\|[g(f)]^2\|_{\wz\ce^{2\az,\,p/2}_{\rho}(\cx)} \le
C\|f\|^2_{\ce^{\az,\,p}_{\rho}(\cx)}$.
\end{thm}

\begin{pf}\rm
By similarity, we only prove the case when $\az>0.$
Let $f\in\ce^{\az,\,p}_{\rho}(\cx)$. By the
homogeneity of $\|\cdot\|_{\ce^{\az,\,p}_{\rho}(\cx)}$
and $\|\cdot\|_{\wz \ce^{\az,\,p}_{\rho}(\cx)}$, we may
assume that $\|f\|_{\ce^{\az,\,p}_{\rho}(\cx)}=1.$ For all balls
$B\equiv B(x_0,\,r)\in{\cd_\rho}$, we need to prove that
\begin{equation}\label{4.6}
\int_B[g(f)(x)]^p\,d\mu(x)\ls[\mu(B)]^{1+\az p}.
\end{equation}
For all $x\in \cx$, write
$$[g(f)(x)]^2=\int_0^{8r} |Q_t(f)(x)|^2\frac {dt}t
+\int_{8r}^\fz |Q_t(f)(x)|^2\frac {dt}t\equiv
[g_{1}(f)(x)]^2+[g_2(f)(x)]^2.$$
By the $L^p(\cx)$-boundedness of
$g$ and \eqref{2.1}, we have
\begin{eqnarray}\label{4.7}
\int_B[g_1(f\chi_{2B})(x)]^p\,d\mu(x)\ls
\int_{2B}|f(x)|^p\,d\mu(x)\ls[\mu(B)]^{1+\az p}.
\end{eqnarray}
By ${\rm (Q)_i}$,
$\gz>\az n$, \eqref{2.1} and the H\"older inequality, we have that
for all $x\in B$ and $t<8r$,
\begin{eqnarray*}
\lf|Q_t\lf(f\chi_{(2B)^\complement}\r)(x)\r|
&&\ls\int_{(2B)^\complement}\frac1{V_t(x)+V(x,\,y)}
\lf(\frac{t}{t+d(x,\,y)}\r)^{\gz}|f(y)|\,d\mu(y)\\
&&\ls\lf(\frac tr\r)^{\gz}\sum_{j=1}^\fz
2^{-j\gz}\frac1{\mu(2^{j+1}B)}
\int_{2^{j+1}B}|f(y)|\,d\mu(y)\\
&&\ls\lf(\frac tr\r)^{\gz}[\mu(B)]^\az \sum_{j=1}^\fz
2^{-j(\gz-\az n )}\ls\lf(\frac tr\r)^{\gz}[\mu(B)]^\az.
\end{eqnarray*}
From this, it follows that
\begin{equation}\label{4.8}
\int_B\lf[g_{1}\lf(f\chi_{(2B)^\complement}\r)(x)\r]^p\,d\mu(x)\ls
\lf(\int_0^{8r}\lf(\frac
tr\r)^{2\gz}\,\frac{dt}t\r)^{p/2}[\mu(B)]^{1+\az p}\ls[\mu(B)]^{1+\az p}.
\end{equation}
 Combining \eqref{4.7} and \eqref{4.8} leads
to that
\begin{equation}\label{4.9}
\int_B[g_1(f)(x)]^p\,d\mu(x)\ls[\mu(B)]^{1+\az p}.
\end{equation}

Applying Lemma 2.1 (ii) and (iii) in \cite{yz08},
we have that for all $x$, $y\in\cx$,
\begin{equation*}
\frac 1{\rho(x)}\gs\frac 1{\rho(y)}
\lf(1+\frac{d(x, y)}{\rho(y)}\r)^{-\frac{k_0}{(1+k_0)}},
\end{equation*}
where $k_0$ is as in Definition \ref{d2.2}.
By this fact, we obtain that for all $x\in B$ and $t\ge 8r$,
$$\frac 1{\rho(x)}\gs \frac 1{\rho(x_0)}
\lf(1+\frac r{\rho(x_0)}\r)^{-\frac{k_0}{(1+k_0)}} \gs \frac
1{\rho(x_0)} \lf(\frac r{\rho(x_0)}\r)^{-\frac{k_0}{(1+k_0)}}.$$
From this, Lemma \ref{l4.1} (i) and \eqref{2.1},  we deduce that
for all $x\in B$,
\begin{eqnarray*}
|Q_t(f)(x)|\ls \lf(\frac{\rho(x)}{t}\r)^{\dz_1}
[\mu(B(x,\,t))]^\az\ls \lf(\frac {\rho(x_0)}t\r)^{\dz_1} \lf(\frac
r{\rho(x_0)}\r)^{\dz_1\frac{k_0}{(1+k_0)}} \lf(\frac{t}{r}\r)^{\az n
}[\mu(B)]^\az,
\end{eqnarray*}
which together with the assumption that $\dz_1>\az n$ implies
that
\begin{eqnarray*}
\int_B[g_2(f)(x)]^p\,d\mu(x) &&\ls [\mu(B)]^{1+\az p}\lf(\frac
r{\rho(x_0)}\r)^{p\dz_1\frac{k_0}{(1+k_0)}}\lf\{\int_{8r}^\fz
\lf(\frac {\rho(x_0)}t\r)^{2\dz_1}\lf(\frac{t}{r}\r)^{2\az n }
\frac {dt}t\r\}^{p/2}\\
&&\ls[\mu(B)]^{1+\az p}\lf(\frac
r{\rho(x_0)}\r)^{p\dz_1\frac{k_0}{(1+k_0)}}
\lf(\frac {\rho(x_0)}r\r)^{p\dz_1}\ls[\mu(B)]^{1+\az p}.
\end{eqnarray*}
This together with \eqref{4.9} gives \eqref{4.6}. Moreover, it
follows from \eqref{4.6} that $g(f)(x)<\fz$ for a.\,e. $x\in\cx$.

Now we assume that $B\equiv B(x_0,\,r)\notin{\cd_\rho}$. We need to prove
that
\begin{equation}\label{4.10}
\int_B\lf\{[g(f)(x)]^2-{\mathop\einf_B}[g(f)]^2\r\}^{p/2}\,d\mu(x)
\ls[\mu(B)]^{1+\az p}.
\end{equation}
To this end, write
\begin{eqnarray*}
[g(f)(x)]^2&&=\int_0^{8r}|Q_t(f)(x)|^2\frac
{dt}t+\int_{8r}^{8\rho(x_0)}\cdots
+\int_{8\rho(x_0)}^\fz\cdots \\
&&\equiv
[g_{r}(f)(x)]^2+[g_{r,\,\rho(x_0)}(f)(x)]^2+[g_{\rho(x_0),\,\fz}(f)(x)]^2.
\end{eqnarray*}
Then \begin{eqnarray*}
&&\int_B\lf\{[g(f)(x)]^2-{\mathop\einf_B}[g(f)]^2\r\}^{p/2}\,d\mu(x)\\
&&\hs\ls\int_B [g_r(f)(x)]^p\,d\mu(x)+
\int_B\lf\{[g_{r,\,\rho(x_0)}(f)(x)]^2-
{\mathop\einf_B}[g_{r,\,\rho(x_0)}(f)]^2\r\}^{p/2}\,d\mu(x)\\
&&\hs\hs\hs+\int_B\lf\{[g_{\rho(x_0),\fz}(f)(x)]^2
-{\mathop\einf_B}[g_{\rho(x_0),\,\fz}(f)]^2\r\}^{p/2}\,d\mu(x)\\
&&\hs\ls\int_B [g_r(f)(x)]^p\,d\mu(x)+\mu(B)\sup_{x,\,y\in
B}\lf|[g_{r,\,\rho(x_0)}(f)(x)]^2
-[g_{r,\,\rho(x_0)}(f)(y)]^2\r|^{p/2} \\
&&\hs\hs\hs+\mu(B)\sup_{x,\,y\in
B}\lf|[g_{\rho(x_0),\,\fz}(f)(x)]^2
-[g_{\rho(x_0),\,\fz}(f)(y)]^2\r|^{p/2}
\equiv {\rm I_1}+{\rm I_2}+{\rm I_3}.
\end{eqnarray*}

Write $f=(f-f_B)\chi_{2B}+(f-f_B)\chi_{(2B)^\complement}+f_B\equiv
f_1+f_2+f_B.$ By the $L^p(\cx)$-boundedness of $g(\cdot)$ and \eqref{2.1},
we have
\begin{equation}\label{4.11}
\dint_B[g_r(f_1)(x)]^p\,d\mu(x)\ls
\int_{2B}|f(x)-f_B|^p\,d\mu(x)\ls[\mu(B)]^{1+\az p}.
\end{equation}
Using ${\rm (Q)_i}$, \eqref{2.1}, the H\"older inequality, Lemma \ref{l2.4} (ii)
and $\gz>\az n$, we obtain that for all $x\in B$,
\begin{eqnarray*}
|Q_t(f_2)(x)|&&\ls\int_{(2B)^\complement}\frac1{V_t(x)+V(x,\,y)}
\lf(\frac{t}{t+d(x,\,y)}\r)^{\gz}
|f(y)-f_B|\,d\mu(y)\\
&&\ls\lf(\frac tr\r)^{\gz}\sum_{j=1}^\fz
2^{-j\gz}\frac1{\mu(2^{j+1}B)} \int_{2^{j+1}B}
[|f(y)-f_{2^{j+1}B}|+|f_{2^{j+1}B}-f_B|]\,d\mu(y)\\
&&\ls\lf(\frac tr\r)^{\gz}[\mu(B)]^\az\sum_{j=1}^\fz 2^{-j(\gz-\az n)}
\ls\lf(\frac tr\r)^{\gz}[\mu(B)]^\az,
\end{eqnarray*}
from which it follows that
\begin{equation}\label{4.12}
\int_B[g_{r}(f_2)(x)]^p\,d\mu(x)\ls [\mu(B)]^{1+\az
p}\lf(\int_0^{8r}\lf(\frac tr\r)^{2\gz}\,\frac{dt}t\r)^{p/2}\ls[\mu(B)]^{1+\az p}.
\end{equation}

Recall that for all $x\in B$, $\rho(x)\sim\rho(x_0)$ (see \eqref{3.6}). By this,
${\rm (Q)_{iii}}$ and Lemma \ref{l2.4} (i), we have that for all $x\in B$,
$$|Q_t(f_B)(x)|\ls \lf(\frac t{t+\rho(x)}\r)^{\dz_2}|f_B|\ls \lf(\frac
t{\rho(x_0)}\r)^{\dz_2}\lf(\frac{\rho(x_0)}{r}\r)^{\az n }
[\mu(B)]^\az.$$ This together with $\dz_2>3\az n$ and
$r<\rho(x_0)$ implies that
\begin{eqnarray*}
\int_B[g_r(f_B)(x)]^p\,d\mu(x)&&\ls [\mu(B)]^{1+\az
p}\lf(\frac{\rho(x_0)}{r}\r)^{\az pn }\lf( \int_0^{8r}\lf(\frac
t{\rho(x_0)}\r)^{2\dz_2}
\,\frac{dt}t\r)^{p/2}\ls [\mu(B)]^{1+\az p}.
\end{eqnarray*}
Combining this, \eqref{4.11} and \eqref{4.12}
yields ${\rm I_1}\ls [\mu(B)]^{1+\az p}$.

Since $\gz>\az n$, by Lemma
\ref{l4.1}, \eqref{2.1} and $\rho(x_0)\sim\rho(x)$ for all $x\in B$, we have that
for all $x$, $y\in B$ and $t\in[8\rho(x_0), \fz)$,
$$|Q_t(f)(x)-Q_t(f)(y)|
\ls\lf(\frac{d(x,\,y)}{t}\r)^{\bz}[\mu(B(x,\,t))]^\az
\ls\lf(\frac{r}{t}\r)^{\bz-\az n }[\mu(B)]^\az,$$ and
$$|Q_t(f)(x)|\ls\lf(\frac{\rho(x_0)}{t}\r)^{\dz_1}[\mu(B(x,\,t))]^\az\ls
\lf(\frac{\rho(x_0)}{t}\r)^{\dz_1}\lf(\frac{t}{r}\r)^{\az n
}[\mu(B)]^\az.$$
By these inequalities and $\bz\ge3\az n$, we see
that for all $x$, $y\in B$,
\begin{eqnarray*}
&&\lf\{[g_{\rho(x_0),\,\fz}(f)(x)]^2-[g_{\rho(x_0),\,\fz}(f)(y)]^2\r\}\\
&&\hs\le\int_{8\rho(x_0)}^\fz
|Q_t(f)(x)+Q_t(f)(y)||Q_t(f)(x)-Q_t(f)(y)|\,\frac{dt}t\\
&&\hs\le\int_{8\rho(x_0)}^\fz
\lf(\frac{\rho(x_0)}{t}\r)^{\dz_1}\lf(\frac{r}{t}\r)^{\bz-2\az n
}[\mu(B)]^{2\az}\,\frac{dt}t\ls [\mu(B)]^{2\az},
\end{eqnarray*}
which implies that ${\rm I_3}\ls [\mu(B)]^{1+\az p}.$

By Lemma \ref{l4.1} (i), \eqref{2.1}
and the fact that for all $x\in B$, $\rho(x_0)\sim\rho(x)$,
we have that for all $t\in [8r, 8\rho(x_0))$ and $x\in B$,
\begin{eqnarray*}
|Q_t(f)(x)|\ls[\mu(B(x,\,t))]^\az\ls\lf(\frac t{r}\r)^{\az
n }[\mu(B)]^\az.
\end{eqnarray*}
Thus the fact that $\bz\ge3\az n$ implies that for all $x$, $y\in B$,
\begin{eqnarray*}
&&\lf\{[g_{r,\,\rho(x_0)}(f)(x)]^2-[g_{r,\,\rho(x_0)}(f)(y)]^2\r\}\\
&&\hs\le\int_{8r}^{8\rho(x_0)}
|Q_t(f)(x)+Q_t(f)(y)||Q_t(f)(x)-Q_t(f)(y)|\,\frac{dt}t\\
&&\hs\ls[\mu(B)]^{\az}
\int_{8r}^{8\rho(x_0)} \lf(\frac t{r}\r)^{\az n}|Q_t(f)(x)-Q_t(f)(y)|\,\frac{dt}t.
\end{eqnarray*}
Let $t\in[8r, 8\rho(x_0))$, $x$, $y\in B$. We write
\begin{eqnarray*}
&&|Q_t(f)(x)-Q_t(f)(y)|\\
&&\hs\le \lf|\dint_\cx\lf[Q_t(x, z)-Q_t(y,
z)\r][f(z)-f_B]\,d\mu(z)\r|
+|f_B|\lf|\dint_\cx\lf[Q_t(x, z)-Q_t(y, z)\r]\,d\mu(z)\r|\\
&&\hs\equiv{\rm H_1}+{\rm H_2}.
\end{eqnarray*}
By ${\rm (Q)_{ii}}$, $t\in [8r, 8\rho(x_0))$, \eqref{2.1}
and Lemma \ref{l2.4} (ii), we see that for all $x\in B$,
\begin{eqnarray*}
{\rm H_1}&&\ls\int_\cx
\lf(\frac{d(x,\,y)}{t+d(x,\,z)}\r)^\bz\frac1{V_t(x)+V(x,\,z)}
\lf(\frac{t}{t+d(x,\,z)}\r)^{\gz}|f(z)-f_B|\,d\mu(z)\\
&&\ls\sum_{j=0}^\fz\frac{r^\bz t^{\gz}}{(t+2^{j-1}r)^{\bz+\gz}}
\frac1{\mu(2^{j+1}B)}\int_{2^{j+1}B}\lf\{|f(z)-f_{2^{j+1}B}|
+|f_{2^{j+1}B}-f_B|\r\}\,d\mu(z)\\
&&\ls\sum_{j=0}^\fz\frac{r^\bz
t^{\gz}}{(t+2^jr)^{\bz+\gz}}2^{j\az n}[\mu(B)]^\az.
\end{eqnarray*}
From this, we deduce that
\begin{equation*}
\dint_{8r}^{8\rho(x_0)}\lf(\frac t{r}\r)^{\az n}{\rm H_1}\,\frac{dt}t \ls [\mu(B)]^\az
\sum_{j=0}^\fz2^{j\az n}\dint_{8r}^{8\rho(x_0)}
\lf(\frac rt\r)^{\bz-\az n}\frac{t^{\gz+\bz-1}}{(t+2^jr)^{\bz+\gz}}\,dt\ls [\mu(B)]^\az.
\end{equation*}

By Lemma \ref{l2.4} (i), ${\rm (Q)_{iii}}$,
$\bz\ge3\az n$, $\dz_2>3\az n$ and the fact that for all $z\in B$,
$\rho(x_0)\sim\rho(z)$,
we have that  for
$\mu$-a.\,e.\,$x$, $y\in B$,
\begin{eqnarray*}
&&\dint_{8r}^{8\rho(x_0)}\lf(\frac t{r}\r)^{\az n}{\rm H_2}\,\frac{dt}t\\
&&\le\dint_{8r}^{8\rho(x_0)}\lf(\frac t{r}\r)^{\az n}
\lf(\frac{\rho(x_0)}r\r)^{\az n}[\mu(B)]^\az
\lf|Q_t(1)(x)-Q_t(1)(y)\r|^{\frac23}
\lf(\frac{t}{\rho(x_0)}\r)^{\frac{\dz_2}3}\,\frac{dt}t\\
&&\ls\dint_{8r}^{8\rho(x_0)}
\lf(\frac{\rho(x_0)}r\r)^{\az
n}[\mu(B)]^\az\lf(\frac r{t}\r)^{\frac{\bz}3}
\lf(\frac{t}{\rho(x_0)}\r)^{
\frac{\dz_2}3}\,\frac{dt}t\\
&&\ls\dint_{8r}^{\rho(x_0)}[\mu(B)]^\az
\lf(\frac {t}{\rho(x_0)}\r)^{\frac{\dz_2}3-\az n}\,\frac{dt}t
\ls[\mu(B)]^\az.
\end{eqnarray*}
This finishes the proof of Theorem \ref{t4.1}.
\end{pf}

As a consequence of Theorem \ref{t4.1}, we have the following
conclusion.

\begin{cor}\label{c4.1}
With the assumptions same as in Theorem \ref{t4.1},
there exists a positive constant $C$ such for
all $f\in\ce^{\az,\,p}_{\rho}(\cx)$,
$g(f)\in\wz\ce^{\az,\,p}_{\rho}(\cx)$ and
$\|g(f)\|_{\wz\ce^{\az,\,p}_{\rho}(\cx)} \le
C\|f\|_{\ce^{\az,\,p}_{\rho}(\cx)}$.
\end{cor}

\begin{pf}\rm
Since
\begin{equation*}
0\le g(f)-{\mathop\einf_Bg(f)}\le\lf\{[g(f)]^2-
{\mathop\einf_B[g(f)]^2}\r\}^{1/2},
\end{equation*}
applying \eqref{4.10}, we have that for all balls $B\notin\cd_\rho$,
\begin{eqnarray}\label{4.13}
&&\lf\{\frac1{[\mu(B)]^{1+\az p}}\int_B \lf[g(f)(x)-
{\mathop\einf_{B}}g(f)\r]^p\,d\mu(x)\r\}^{1/p}\\
&&\hs\ls \lf\{\frac1{[\mu(B)]^{1+\az p}}\dint_B \lf\{\lf[g(f)(x)\r]^2-
{\mathop\einf_B}[g(f)]^2\r\}^{p/2}\,d\mu(x)\r\}^{1/p}
\ls\|f\|_{\ce^{\az,\,p}_{\rho}(\cx)}.\nonumber
\end{eqnarray}
On the other hand, by \eqref{4.6},
we obtain that for all balls $B\in\cd_\rho$,
\begin{equation*}
\lf\{\frac1{[\mu(B)]^{1+\az p}}\int_B \lf[g(f)(x)
\r]^p\,d\mu(x)\r\}^{1/p}\ls\|f\|_{\ce^{\az,\,p}_{\rho}(\cx)},
\end{equation*}
which together with \eqref{4.13} completes the proof of Corollary
\ref{c4.1}.
\end{pf}

\begin{rem}\label{r4.1}\rm
(i) If $\az=0$ and $\cx$ is an RD-space, Theorem \ref{t4.1} and
Corollary \ref{c4.1} were already obtained in \cite{yyz}.

(ii) We point out that Remark \ref{r3.1} (i) is also suitable
to Theorem \ref{t4.1} and Corollary \ref{c4.1}.
\end{rem}

\section{Several applications}\label{s5}

\hskip\parindent This section is divided into
Subsections \ref{s5.1} through \ref{s5.4}.
We apply the results obtained in Sections \ref{s3} and \ref{s4}, respectively,
to the Schr\"odinger operator or the
degenerate Schr\"odinger operator on $\rd$, the
sub-Laplace Schr\"odinger operator on Heisenberg groups or
on connected and simply connected nilpotent Lie groups.

\subsection{Schr\"odinger operators on $\rd$}\label{s5.1}

\hskip\parindent Let $d\ge3$ and $\rd$ be the
$d$-dimensional Euclidean space endowed with the Euclidean norm
$|\cdot|$ and the Lebesgue measure $dx$. Denote the Laplacian
$\sum_{j=1}^d\frac{\partial^2}{\partial x_j^2}$ on $\rd$ by
$\Delta$ and the corresponding heat (Gauss) semigroup
$\{e^{t\Delta}\}_{t>0}$ by $\{\wz T_t\}_{t>0}$. Let $V$ be a
nonnegative locally integrable function on $\rd$,
$\cl\equiv-\Delta+V$ the Schr\"odinger operator and
$\{T_t\}_{t>0}$ the corresponding semigroup with
kernels $\{T_t(x,\,y)\}_{t>0}$.
Moreover, for all $t>0$ and $x,\,y\in\rd$, set
$$Q_t(x,\,y)\equiv t^2\frac{d T_s(x,\,y)}{ds}\Bigg|_{s=t^2}.$$
Let $q \in(d/2, d]$, $V\in\cb_q(\rd, |\cdot|, dx)$ and $\rho$ be as
in \eqref{2.3}. Then we have the following estimates; see
\cite{d05,dz03,dgmtz05}.

\begin{prop}\label{p5.1}
Let $q \in(d/2, d]$, $\bz\in(0, 2-d/q)$ and $N\in\nn$.
Then there exist  positive constants $\wz C$ and $C$,
where $C$ is independent of $N$, such that
for all $t\in(0,\,\fz)$ and $x,\,x',\,y\in\cx$ with $d(x,\,x')\le
\sqrt t/2$,
\begin{enumerate}
\vspace{-0.2cm}
\item[(i)] $|T_t(x,\,y)|\le \wz Ct^{-d/2}\exp\{-\frac
{|x-y|^2}{Ct}\}[\frac {\rho(x)}{\sqrt t+\rho(x)}]^N
[\frac {\rho(y)}{\sqrt t+\rho(y)}]^N$,
\vspace{-0.2cm}
\item[(ii)]  $|T_t(x,\,y)-T_t(x',\,y)|\le\wz C
 [\frac {|x-x'|}{\sqrt t}]^{\bz}t^{-d/2}\exp \{-\frac
{|x-y|^2}{Ct}\}[\frac {\rho(x)}{\sqrt t+\rho(x)}]^N
[\frac {\rho(y)}{\sqrt t+\rho(y)}]^N$,
\vspace{-0.2cm}
\item[(iii)]  $|T_t(x,\,y)-\wz T_t(x,\,y)|\le\wz C
 [\frac {\sqrt t}{\sqrt t+\rho(x)}]^{2-d/q}t^{-d/2}\exp\{-\frac
{|x-y|^2}{Ct}\}$;
\end{enumerate}
\vspace{-0.1cm}
\noindent and for all $t\in(0,\,\fz)$ and $x,\,x',\,y\in\cx$ with $d(x,\,x')\le
t/2$,
\begin{enumerate}
\vspace{-0.3cm}
\item[(iv)] $|Q_t(x,\,y)|\le  \wz C t^{-d}\exp \{-\frac
{|x-y|^2}{Ct^2}\}[\frac {\rho(x)}{t+\rho(x)}]^N[\frac {\rho(y)}{t+\rho(y)}]^N$,
\vspace{-0.3cm}
\item[(v)]  $|Q_t(x,\,y)-Q_t(x',\,y)|\le\wz C
 [\frac {|x-x'|}{t}]^{\bz}t^{-d}\exp \{-\frac
{|x-y|^2}{Ct^2}\}[\frac {\rho(x)}{t+\rho(x)}]^N
[\frac {\rho(y)}{t+\rho(y)}]^N$,
\vspace{-0.3cm}
\item[(vi)] $|\int_\rd Q_t(x,\,y)d\mu(y)|\le \wz C
[\frac t{\rho(x)}]^{2-d/q}[\frac{\rho(x)}{t+\rho(x)}]^N$.
\end{enumerate}
\end{prop}

Let $q_1$, $q_2\in(d/2,\fz]$ with $q_1<q_2$.
Observe that $\cb_{q_2}(\rd)\subset \cb_{q_1}(\rd)$.
Therefore, Proposition \ref{p5.1} holds for all
$q\in(d/2, \fz]$. On the other hand,
recall that $\{\wz T_{t^2}\}_{t>0}$ satisfies that
for all $t\in(0, \fz)$, $\wz T_{t^2}(1)=1$
(see \cite{d05,dz03}). Thus $\{T_{t^2}\}_{t>0}$
satisfies the assumptions \eqref{3.1} through \eqref{3.3}.
Moreover, the $L^2(\rd)$-boundedness of $g$-function $g(\cdot)$ was obtained in
\cite{d05}. Using this, (iv) and (v) of Proposition \ref{p5.1} and
the vector-valued Calder\'on-Zygmund theory (see, for example,
\cite{s93}), we obtain the $L^p(\rd)$-boundedness of $g(\cdot)$ for
$p\in(1,\,\fz)$.
Then by applying this fact and Proposition \ref{p5.1},
Theorems \ref{t3.1}, \ref{t3.2}, \ref{t4.1}
and Corollary \ref{c4.1}, we have the following result.

\begin{prop}\label{p5.2}
Let $q\in (d/2,\fz]$, $p\in(1,\,\fz)$, $V\in\cb_q(\rd, |\cdot|, dx)$
and $\rho$ be as in \eqref{2.3}.
\begin{enumerate}
\item[(i)] If $\az\in(-\fz,\,1/d-1/(2q))$,
then there exists a positive constant $C$ such that for all
$f\in\ce^{\az,\,p}_{\rho}(\rd)$, $T^+(f)$, $P^+(f)\in \wz\ce^{\az,\,p}_{\rho}(\rd)$ and
$$\|T^+(f)\|_{\wz\ce^{\az,\,p}_{\rho}(\rd)}
+\|P^+(f)\|_{\wz\ce^{\az,\,p}_{\rho}(\rd)}\le
C\|f\|_{\ce^{\az,\,p}_{\rho}(\rd)}.$$

\item[(ii)] If $\az\in(-\fz,\,2/(3d)-1/(3q))$,
then there exists a positive constant $C$ such that for all
$f\in\ce^{\az,\,p}_{\rho}(\rd)$,
$[g(f)]^2\in\wz\ce^{2\az,\,p/2}_{\rho}(\rd)$ with
$\|[g(f)]^2\|_{\wz\ce^{2\az,\,p/2}_{\rho}(\rd)} \le
C\|f\|^2_{\ce^{\az,\,p}_{\rho}(\rd)}$, and
$g(f)\in\wz\ce^{\az,\,p}_{\rho}(\rd)$ with
$\|g(f)\|_{\wz\ce^{\az,\,p}_{\rho}(\rd)} \le
C\|f\|_{\ce^{\az,\,p}_{\rho}(\rd)}$.
\end{enumerate}
\end{prop}

\subsection{Degenerate Schr\"odinger operators on  $\rd$}\label{s5.2}

\hskip\parindent Let $d\ge 3$ and $\rd$ be the $d$-dimensional Euclidean space
endowed with the Euclidean norm $|\cdot|$ and the Lebesgue measure $dx$.
Recall that a nonnegative locally integrable function $w$ is said to
be an $A_2(\rd)$ weight in the sense of Muckenhoupt if
$$\sup_{B\subset \rd}\lf\{\frac1{|B|}\int_B w(x)\,dx\r\}^{1/2}
\lf\{\frac1{|B|}\int_B [w(x)]^{-1}\,dx\r\}^{1/2}<\fz,$$
where the supremum is taken over all the balls in $\rd$.
Observe that if we set $w(E)\equiv\int_Ew(x)dx$ for any measurable set $E$, then
 there exist positive constants $C,\,Q$ and $\kz$ such that for all $x\in\rd$, $\lz>1$ and $r>0$,
$$C^{-1}\lz^\kz w(B(x,\,r))\le w(B(x,\,\lz r))\le C\lz^Qw(B(x,\,r)),$$
namely, the measure $w(x)\,dx$ satisfies \eqref{2.1}.
Thus $(\rd,\,|\cdot|,\,w(x)\,dx)$ is a space of homogeneous type.

Let $w\in A_2(\rd)$ and
$\{a_{i,\,j}\}_{1\le i,\,j\le d}$ be a real symmetric matrix function satisfying that
for all $x,\,\xi\in\rd$,
$$C^{-1}|\xi|^{2}\le\sum_{1\le i,\,j\le d}a_{i,\,j}(x)\xi_i\overline \xi_j
\le C|\xi|^2.$$ Then the degenerate elliptic operator
$\cl_0$ is defined by
$$\cl_0 f(x)\equiv-\frac1{w(x)}\sum_{1\le i,\,j\le d}
\partial_i(a_{i,\,j}(\cdot)\partial_j f)(x),$$
where $x\in\rd$. Denote by $\{\wz T_t\}_{t>0}\equiv \{e^{-t\cl_0}\}_{t>0}$
the semigroup generated by $\cl_0$.

Let $V$ be a nonnegative locally integrable function on $w(x)\,dx$.
Define the  degenerate Schr\"odinger operator by $\cl \equiv\cl_0+V.$
Then $\cl$ generates a semigroup $\{T_t\}_{t>0}\equiv\{e^{-t\cl}\}_{t>0}$ with
kernels $\{T_t(x,\,y)\}_{t>0}$.
Moreover, for all $t\in(0, \fz)$ and $x,\,y\in\rd$, set
$$Q_t(x,\,y)\equiv t^2\frac{d T_s(x,\,y)}{ds}\Bigg|_{s=t^2}.$$
Let $q \in(Q/2, Q]$,
$V\in\cb_q(\rd, |\cdot|, w(x)\,dx)$ and $\rho$ be as in \eqref{2.3}.
Then $\{T_t(\cdot, \cdot)\}_{t>0}$ and
$\{Q_t(\cdot, \cdot)\}_{t>0}$ satisfy Proposition \ref{p5.1}
with $t^{-d/2}$ replaced by $[V_{\sqrt t}(x)]^{-1}$,
$t^{-d}$ by $[V_{ t}(x)]^{-1}$ and $d $ by $Q$.
In fact, the corresponding Proposition \ref{p5.1}
(i) and (iii) here were given in \cite{d05}.
The proof of the corresponding Proposition \ref{p5.1} (ii)
here is similar to that of Proposition \ref{p5.1}; see \cite{dz03}.
The proofs of the corresponding Proposition \ref{p5.1}
(iv), (v) and (vi) here are similar to that of Proposition 4 of \cite{dgmtz05}.
We omit the details here.

Recall that $\{\wz T_{t^2}\}_{t>0}$ satisfies that
for all $t\in(0, \fz)$, $\wz T_{t^2}(1)=1$;
see, for example, \cite{hs01}. Thus $\{T_{t^2}\}_{t>0}$
satisfies the assumptions \eqref{3.1} through \eqref{3.3}.
Moreover, the $L^2(\rd)$-boundedness of
$g(\cdot)$ can be obtained by the same argument
as in Lemma 3 of \cite{d05}. Using this, (iv) and (v) of Proposition \ref{p5.1} and
the vector-valued Calder\'on-Zygmund theory,
we obtain the $L^p(\rd)$-boundedness of $g(\cdot)$ for
$p\in(1,\,\fz)$.
Then by applying these facts and
Theorems \ref{t3.1}, \ref{t3.2}, \ref{t4.1} and
Corollary \ref{c4.1}, we have the following result.

\begin{prop}\label{p5.3}
Let $w\in A_2(\rd)$, $q\in (Q/2,\fz]$, $p\in(1,\,\fz)$,
$V\in \cb_q(\rd,\,|\cdot|,\,w(x)\,dx)$
and $\rho$ be as in \eqref{2.3} with $d\mu=w(x)\,dx$.
\begin{enumerate}
\item[(i)] If  $\az\in(-\fz,\,1/Q-1/(2q))$,
then there exists a positive constant $C$ such that for all
$f\in\ce^{\az,\,p}_{\rho}(w(x)\,dx)$, $T^+(f)$,
$P^+(f)\in\wz\ce^{\az,\,p}_{\rho}(w(x)\,dx)$ and
$$\|T^+(f)\|_{\wz\ce^{\az,\,p}_{\rho}(w(x)\,dx)}
+\|P^+(f)\|_{\wz\ce^{\az,\,p}_{\rho}(w(x)\,dx)}\le
C\|f\|_{\ce^{\az,\,p}_{\rho}(w(x)\,dx)}.$$

\item[(ii)] If $\az\in(-\fz,\,2/(3Q)-1/(3q))$,
then there exists a positive constant $C$ such that for all
$f\in\ce^{\az,\,p}_{\rho}(w(x)\,dx)$,
$[g(f)]^2\in\wz\ce^{2\az,\,p/2}_{\rho}(w(x)\,dx)$ with
$\|[g(f)]^2\|_{\wz\ce^{2\az,\,p/2}_{\rho}(w(x)\,dx)} \le
C\|f\|^2_{\ce^{\az,\,p}_{\rho}(w(x)\,dx)}$, and
$g(f)\in\wz\ce^{\az,\,p}_{\rho}(w(x)\,dx)$ with
$$\|g(f)\|_{\wz\ce^{\az,\,p}_{\rho}(w(x)\,dx)} \le
C\|f\|_{\ce^{\az,\,p}_{\rho}(w(x)\,dx)}.$$
\end{enumerate}
\end{prop}

\subsection{Schr\"odinger operators on Heisenberg groups}\label{s5.3}

\hskip\parindent The $(2n+1)$-dimensional Heisenberg group $\hh^n$ is
a connected and simply connected nilpotent Lie group
with the underlying manifold $\rr^{2n}\times \rr$ and the multiplication
$$(x,\,s)(y,\,s)=\lf(x+y,\,t+s+2\sum_{j=1}^n[x_{n+j}y_j-x_jy_{n+j}]\r).$$
The homogeneous norm on $\hh^n$ is defined by
$|(x,\,t)|=(|x|^4+|t|^2)^{1/4}$ for all $(x,\,t)\in\hh^n$, which
induces a left-invariant metric
$d((x,\,t),\,(y,\,s))=|(-x,\,-t)(y,\,s)|$. Moreover, there exists a
positive constant $C$ such that $|B((x,\,t),\,r)|=Cr^Q,$ where
$Q=2n+2$ is the homogeneous dimension of $\hh^n$ and
$|B((x,\,t),\,r)|$ is the Lebesgue measure of the ball
$B((x,\,t),\,r)$. The triplet $(\hh^n,\, d,\,dx)$ is a space of homogeneous type.

A basis for the Lie algebra of left invariant vector fields on
$\hh^n$ is given by
$$
X_{2n+1}=\frac{\partial}{\partial t}, \hs
X_j=\frac{\partial}{\partial x_j}+2x_{n+j}\frac{\partial}{\partial
t}, \hs X_{n+j}= \frac{\partial}{\partial
x_{n+j}}-2x_j\frac{\partial}{\partial t},\ j=1,\,\cdots,\,n.$$ All
non-trivial commutators are $[X_j,\,X_{n+j}]=-4X_{2n+1}$,
$j=1,\,\cdots,\,n$. The sub-Laplacian has the form
$\Delta_{\hh^n}=\sum_{j=1}^{2n}X_j^2.$

Let $V$ be a nonnegative locally integrable function on $\hh^n$.
Define the sub-Laplacian Schr\"odinger operator by
$\cl\equiv-\Delta_{\hh^n}+V.$ Denote by $\{T_t\}_{t>0}\equiv
\{e^{-t\cl}\}_{t>0}$ with
kernels $\{T_t(x,\,y)\}_{t>0}$ and by $\{\wz T_t\}_{t>0}\equiv
\{e^{t\Delta_{\hh^n}}\}_{t>0}$.
Moreover, for all $t\in(0, \fz)$ and $x,\,y\in\rd$, set
$$Q_t(x,\,y)\equiv t^2\frac{d T_s(x,\,y)}{ds}\Bigg|_{s=t^2}.$$

Let $V\in\cb_q({\hh^n}, d, dx)$ with $q\in(n+1, 2n+2]$ and $\rho$ be as in \eqref{2.3}.
Then $\{T_t(\cdot, \cdot)\}_{t>0}$ and $\{Q_t(\cdot, \cdot)\}_{t>0}$ satisfy Proposition \ref{p5.1}
with $d$ replaced by $2(n+2)$ and $|x-y|$ replaced by $d(x,\,y)$; see \cite{ll08}.

Observe that $\{\wz T_{t^2}\}_{t>0}$ satisfies that for all $t\in(0, \fz)$,
$\wz T_{t^2}(1)=1$ (see also \cite{yyz}). Thus $\{T_{t^2}\}_{t>0}$
satisfies the assumptions \eqref{3.1} through \eqref{3.3}.
Moreover, the $L^2(\hh^n)$-boundedness of
$g(\cdot)$ was obtained in \cite{ll08}.
Using this, (iv) and (v) of Proposition \ref{p5.1} and
the vector-valued Calder\'on-Zygmund theory,
we obtain the $L^p(\hh^n)$-boundedness of $g(\cdot)$ for
$p\in(1,\,\fz)$.
Then by applying these facts and
Theorems \ref{t3.1}, \ref{t3.2}, \ref{t4.1} and Corollary \ref{c4.1},  we have the following
conclusions.

\begin{prop}\label{p5.4}
Let $q\in (n+1,\fz]$, $p\in(1,\,\fz)$, $V\in\cb_q({\hh^n}, d, dx)$
and $\rho$ be as in \eqref{2.3}.
\begin{enumerate}
\item[(i)] If  $\az\in(-\fz,\,1/(2n+2)-1/(2q))$,
then there exists a positive constant $C$ such that for all
$f\in\ce^{\az,\,p}_{\rho}({\hh^n})$, $T^+(f)$, $P^+(f)\in
\wz\ce^{\az,\,p}_{\rho}({\hh^n})$ and
$$\|T^+(f)\|_{\wz\ce^{\az,\,p}_{\rho}({\hh^n})}
+\|P^+(f)\|_{\wz\ce^{\az,\,p}_{\rho}({\hh^n})}\le
C\|f\|_{\ce^{\az,\,p}_{\rho}({\hh^n})}.$$

\item[(ii)] If $\az\in(-\fz,\,1/(3n+3)-1/(3q))$,
then there exists a positive constant $C$ such that for all
$f\in\ce^{\az,\,p}_{\rho}({\hh^n})$,
$[g(f)]^2\in\wz\ce^{2\az,\,p/2}_{\rho}({\hh^n})$ with
$\|[g(f)]^2\|_{\wz\ce^{2\az,\,p/2}_{\rho}({\hh^n})} \le
C\|f\|^2_{\ce^{\az,\,p}_{\rho}({\hh^n})},$ and
$g(f)\in\wz\ce^{\az,\,p}_{\rho}({\hh^n})$ with
$\|g(f)\|_{\wz\ce^{\az,\,p}_{\rho}({\hh^n})} \le
C\|f\|_{\ce^{\az,\,p}_{\rho}({\hh^n})}$ .
\end{enumerate}
\end{prop}

\subsection{Schr\"odinger operators on connected and
simply connected nilpotent Lie groups}\label{s5.4}

\hskip\parindent Let $\mathbb G$ be a connected and simply
connected nilpotent Lie group and
$X\equiv\{X_1,\,\cdots,\,X_k\}$ left invariant vector fields on
$\bbg$ satisfying the H\"ormander condition that
$\{X_1,\,\cdots,\,X_k\}$ together with their commutators of order
$\le m$ generates the tangent space of $\bbg$ at each point of
$\bbg$. Let $d$ be the Carnot-Carath\'eodory (control) distance on
$\bbg$ associated to $\{X_1,\,\cdots,\,X_k\}$. Fix a left invariant
Haar measure $\mu$ on $\bbg$. Then for all $x\in \bbg$,
$V_r(x)=V_r(e)$; moreover, there exist $\kz$, $D\in(0, \fz)$ with
$\kz\le D$ such that for all $x\in\bbg$,
$C^{-1} r^\kz\le V_r(x)\le Cr^\kz$ when
$r\in(0, 1]$, and $C^{-1} r^D\le V_r(x)\le Cr^D$ when $r\in(1, \fz)$; see
\cite{nsw85} and \cite{v88}. Thus $({\mathbb
\bbg},\,d,\,\mu)$ is a space of homogeneous type.

The sub-Laplacian is given by $\Delta_\bbg\equiv\sum_{j=1}^kX_j^2.$
Denote by $\{\wz T_t\}_{t>0}\equiv\{e^{t\Delta_\bbg}\}_{t>0}$ the
semigroup generated by $-\Delta_\bbg$.

Let $V$ be a nonnegative locally integrable function on $\bbg$. Then
the sub-Laplace Schr\"odinger operator $\cl$ is defined by
$\cl\equiv-\Delta_\bbg+V.$ The operator  $\cl$ generates a semigroup
of operators $\{T_t\}_{t>0}\equiv\{e^{-t\cl}\}_{t>0}$, whose kernels
are denoted by $\{T_t(x,\,y)\}_{t>0}$.
Define the radial maximal operator $T^+$ by
$T^+(f)(x)\equiv\sup_{t>0}|e^{-t\cl}(f)(x)|$ for all $x\in\bbg$.

Let $q>D/2$,  $V\in \cb_q(\bbg,\,d,\,\mu)$ and
$\rho$ be as in \eqref{2.3}.
For all $x,\,y\in\bbg$ and $t\in(0, \fz)$, define
$$Q_t(x,\,y)\equiv t^2\frac{d}{ds}\Big|_{s=t^2}T_s(x,y).$$
Then $\{T_t(\cdot, \cdot)\}_{t>0}$ and
$\{Q_t(\cdot, \cdot)\}_{t>0}$ satisfy Proposition \ref{p5.1}
with $t^{-d}$ replaced by $[V_{ t}(x)]^{-1}$,
$t^{-d/2}$ by $[V_{\sqrt t}(x)]^{-1}$ and $d $ by $D$ (see \cite{yz08,yyz}).
Observe that $\{\wz T_{t^2}\}_{t>0}$ satisfies that for all $t\in(0, \fz)$,
$\wz T_{t^2}(1)=1$; see, for example, \cite{v88}. Thus $\{T_{t^2}\}_{t>0}$
satisfies the assumptions \eqref{3.1} through \eqref{3.3}.
Moreover, the $L^2(\bbg)$-boundedness of
$g(\cdot)$ can be obtained
by the same argument as in Lemma 3 in \cite{d05}.
Using this, (iv) and (v) of Proposition \ref{p5.4} and
the vector-valued Calder\'on-Zygmund theory,
we obtain the $L^p(\bbg)$-boundedness of $g(\cdot)$ for
$p\in(1,\,\fz)$.
Then by applying these facts and
Theorems \ref{t3.1}, \ref{t3.2}, \ref{t4.1}
and Corollary \ref{c4.1},  we have the following
conclusions.

\begin{prop}\label{p5.5}
Let $q\in (D/2, \fz]$, $p\in(1,\,\fz)$, $V\in \cb_q(\bbg,\,d,\,\mu)$ and
$\rho$ be as in \eqref{2.3}.
\begin{enumerate}
\item[(i)] If $\az\in(-\fz,\,1/D-1/(2q))$,
then there exists a positive constant $C$ such that for all
$f\in\ce^{\az,\,p}_{\rho}(\bbg)$, $T^+(f)$, $P^+(f)\in\wz\ce^{\az,\,p}_{\rho}(\bbg)$
and
$$\|T^+(f)\|_{\wz\ce^{\az,\,p}_{\rho}(\bbg)}
+\|P^+(f)\|_{\wz\ce^{\az,\,p}_{\rho}(\bbg)}\le
C\|f\|_{\ce^{\az,\,p}_{\rho}(\bbg)}.$$

\item[(ii)] If $\az\in(-\fz,\,2/(3D)-1/(3q))$,
then there exists a positive constant $C$ such that for all
$f\in\ce^{\az,\,p}_{\rho}(\bbg)$,
$[g(f)]^2\in\wz\ce^{2\az,\,p/2}_{\rho}(\bbg)$ with
$\|[g(f)]^2\|_{\wz\ce^{2\az,\,p/2}_{\rho}(\bbg)} \le
C\|f\|^2_{\ce^{\az,\,p}_{\rho}(\bbg)}$, and
$g(f)\in\wz\ce^{\az,\,p}_{\rho}(\bbg)$ with
$\|g(f)\|_{\wz\ce^{\az,\,p}_{\rho}(\bbg)} \le
C\|f\|_{\ce^{\az,\,p}_{\rho}(\bbg)}$ .
\end{enumerate}
\end{prop}

\bigskip

\noindent Dachun Yang, Dongyong Yang and Yuan Zhou

\medskip

\noindent School of Mathematical Sciences,
 Beijing Normal University, Laboratory of Mathematics and Complex systems,
Ministry of Education, Beijing 100875, People's Republic of China

\medskip

\noindent{\it E-mail addresses}: \texttt{dcyang@bnu.edu.cn},
\texttt{dyyang@mail.bnu.edu.cn} and

\hspace{2.56cm}\texttt{yuanzhou@mail.bnu.edu.cn}

\section*{References}

\begin{enumerate}

\bibitem[1]{c63} S. Campanato, Propriet\`a di h\"olderianit\`a
di alcune classi di funzioni, Ann. Scuola Norm. Sup. Pisa (3) 17
(1963), 175-188.

\vspace{-0.28cm}
\bibitem[2]{cr80} R. R. Coifman and R. Rochberg,
Another characterization of BMO, Proc. Amer. Math. Soc. 79 (1980),
249-254.

\vspace{-0.28cm}
\bibitem[3]{cw71} R. R. Coifman and G. Weiss,
 Analyse Harmonique Non-commutative sur Certains Espaces Homog\`enes,
Lecture Notes in Math. 242, Springer, Berlin, 1971.

\vspace{-0.28cm}
\bibitem[4]{cw77}  R. R. Coifman and G. Weiss, Extensions of
Hardy spaces and their use in analysis, Bull. Amer. Math. Soc. 83
(1977), 569-645.

\vspace{-0.3cm}
\bibitem[5]{dxy} X. T. Duong, J. Xiao and L. Yan, Old and new Morrey
spaces with heat kernel bounds, J. Fourier Anal. Appl. 13 (2007),
87-111.

\vspace{-0.28cm}
\bibitem[6]{dz99}
J. Dziuba\'nski and J. Zienkiewicz, Hardy space $H\sp 1$ associated
to Schr\"odinger operator with potential satisfying reverse H\"older
inequality, Rev. Mat. Ibero. 15 (1999), 279-296.


\vspace{-0.28cm}
\bibitem[7]{dz03}
J. Dziuba\'nski and J. Zienkiewicz, $H\sp p$ spaces associated with
Schr\"odinger operators with potentials from reverse H\"older
classes, Colloq. Math. 98 (2003), 5-38.

\vspace{-0.28cm}
\bibitem[8]{d05}
J. Dziuba\'nski, Note on $H\sp 1$ spaces related to degenerate
Schr\"odinger operators, Illinois J. Math. 49 (2005), 1271-1297.

\vspace{-0.28cm}
\bibitem[9]{dgmtz05}
J. Dziuba\'nski, G. Garrig\'os, T. Mart\'inez, J. L. Torrea and J.
Zienkiewicz, $BMO$ spaces related to Schr\"odinger operators with
potentials satisfying a reverse H\"older inequality, Math. Z. 249
(2005), 329-356.

\vspace{-0.28cm}
\bibitem[10]{f83} C. Fefferman, The uncertainty principle,
Bull. Amer. Math. Soc. (N. S.) 9 (1983), 129-206.

\vspace{-0.28cm}
\bibitem[11]{g79}
D. Goldberg, A local version of real Hardy spaces, Duke Math. J. 46
(1979), 27-42.

\vspace{-0.28cm}
\bibitem[12]{hmy2} Y. Han, D. M\"uller and D. Yang,
A theory of Besov and Triebel-Lizorkin spaces on metric measure
spaces modeled on Carnot-Carath\'eodory spaces, Abstr. Appl. Anal.
2008, Art. ID 893409, 250 pp.

\vspace{-0.28cm}
\bibitem[13]{hs01}
W. Hebisch and L. Saloff-Coste, On the relation between elliptic and
parabolic Harnack inequalities, Ann. Inst. Fourier (Grenoble) 51
(2001), 1437-1481.

\vspace{-0.28cm}
\bibitem[14]{hmy07} G. Hu, Y. Meng and D. Yang, Estimates for
Marcinkiewicz integrals in {\rm BMO} and Campanato spaces, Glasg.
Math. J. 49 (2007), 167-187.

\vspace{-0.28cm}
\bibitem[15]{hyy} G. Hu, Da. Yang and Do. Yang,
 $h^1$, $\mathop\mathrm{bmo}$, $\mathop\mathrm{blo}$
and Littlewood-Paley $g$-functions with non-doubling measures, Rev.
Mat. Ibero. 25 (2009), 595-667.

\vspace{-0.28cm}
\bibitem[16]{hl}  J. Huang and H. Liu,
Area integrals associated to Schr\"odinger operators, Submitted.

\vspace{-0.28cm}
\bibitem[17]{n06} E. Nakai, The Campanato, Morrey and H\"older spaces
on spaces of homogeneous type, Studia Math. 176 (2006), 1-19.

\vspace{-0.28cm}
\bibitem[18]{n08} E. Nakai, Orlicz-Morrey spaces and the
Hardy-Littlewood maximal function, Studia Math. 188 (2008), 193-221.

\vspace{-0.28cm}
\bibitem[19]{l07} P. G. Lemari\'e-Rieusset, The Navier-Stokes equations
in the critical Morrey-Campa-nato space, Rev. Mat. Ibero. 23 (2007), 897-930.

\vspace{-0.28cm}
\bibitem[20]{l99}
H. Li, Estimations $L\sp p$ des op\'erateurs de Schr\"odinger sur
les groupes nilpotents, J. Funct. Anal. 161 (1999), 152-218.

\vspace{-0.28cm}
\bibitem[21]{ll08} C. Lin and H. Liu, The BMO-type space $\bmo_\cl$
associated with Schr\"odinger operators on the Heisenberg group,
Preprint.

\vspace{-0.28cm}
\bibitem[22]{ms791} R. A. Mac{\'\i}as and C. Segovia,
Lipschitz functions on spaces of homogeneous type, Adv. in Math. 33
(1979), 257-270.

\vspace{-0.28cm}
\bibitem[23]{nsw85} A. Nagel, E. M. Stein and S. Wainger,
Balls and metrics defined by vector fields I. Basic properties, Acta
Math. 155 (1985), 103-147.

\vspace{-0.28cm}
\bibitem[24]{s95}
Z. Shen, $L\sp p$ estimates for Schr\"odinger operators with certain
potentials, Ann. Inst. Fourier (Grenoble) 45 (1995), 513-546.

\vspace{-0.28cm}
\bibitem[25]{s93} E. M. Stein, Harmonic Analysis:
Real-variable Methods, Orthogonality, and Oscillatory Integrals,
Princeton University Press, Princeton, N. J., 1993.

\vspace{-0.3cm}
\bibitem[26]{st89}
J.-O. Str\"omberg and A. Torchinsky, Weighted Hardy Spaces, Lecture
Notes in Mathematics, 1381, Springer-Verlag, Berlin, 1989.

\vspace{-0.3cm}
\bibitem[27]{p69}J. Peetre, On the theory of ${\mathcal L}_{p,\,\lz}$ spaces,
 J. Funct. Anal. 4 (1969), 71-87.

\vspace{-0.3cm}
\bibitem[28]{tw}
M. H. Taibleson and G. Weiss, The molecular characterization of
certain Hardy spaces. Representation theorems for Hardy spaces,
 pp. 67-149, Ast\'erisque, 77, Soc. Math. France, Paris, 1980.

\vspace{-0.3cm}
\bibitem[29]{t92} H. Triebel, Theory of Function Spaces.
 II, Birkh\"auser Verlag, Basel, 1992.

\vspace{-0.28cm}
\bibitem[30]{v88} N. Th. Varopoulos, Analysis on Lie groups,
J. Funct. Anal. 76 (1988), 346-410.

\vspace{-0.28cm}
\bibitem[31]{v92} N. Th. Varopoulos, L. Saloff-Coste and T. Coulhon,
 Analysis and Geometry on Groups, Cambridge University Press,  Cambridge, 1992.

\vspace{-0.28cm}
\bibitem[32]{yyz} Da. Yang, Do. Yang and Y. Zhou,
Localized BMO and BLO spaces on ${\rm RD}$-spaces
and applications to Schr\"odinger operators, arXiv: 0903.4536.

\vspace{-0.28cm}
\bibitem[33]{yz08} D. Yang and Y. Zhou, Localized Hardy
spaces $H^1$ related to admissible functions on RD-spaces and
applications to Schr\"odinger operators, Trans. Amer. Math. Soc. (to appear).

\vspace{-0.28cm}
\bibitem[34]{z99} J. Zhong, The Sobolev estimates for some
Schr\"odinger type operators, Math. Sci. Res. Hot-Line 3:8 (1999),
1-48.

\end{enumerate}
\end{document}